\documentclass[11pt,a4paper]{article}

\usepackage[T1]{fontenc}
\usepackage{graphicx}
\usepackage{amssymb , amsmath, amsthm}
\usepackage{dsfont}
\usepackage{a4wide}
\usepackage{slashbox}

\newtheorem{thm}{Theorem}[section]

\newtheorem{cor}[thm]{Corollary}
\newtheorem{dfn}{Definition}[section]

\newtheorem{exms}{Examples}[section]

\def\Ker{\mathop{\rm Ker}\nolimits}

\title{Transductive versions of the LASSO\\and the Dantzig Selector}
\author{Pierre Alquier\footnote{LPMA (Univ. Paris 7), CREST-LS.} and Mohamed Hebiri\footnote{LPMA
(Univ. Paris 7), ETH-Z\"urich}}
\date{}

\allowdisplaybreaks

\begin{document}

\maketitle
\begin{abstract}
Transductive methods are useful in prediction problems when the training dataset is composed of a large number of unlabeled observations and a smaller number of labeled observations.
In this paper, we propose an approach for developing transductive prediction procedures that are able to take advantage of the sparsity in the high dimensional linear regression.
More precisely, we define transductive versions of the LASSO \cite{Tibshirani-LASSO} and the Dantzig Selector \cite{Dantzig}.
These procedures combine labeled and unlabeled observations of the training dataset to produce a prediction for the unlabeled observations.
We propose an experimental study of the transductive estimators, that shows that they improve the LASSO and Dantzig Selector in many situations, and particularly in high dimensional problems when the predictors are correlated.
We then provide non-asymptotic theoretical guarantees for these estimation methods.
Interestingly, our theoretical results show that the Transductive LASSO and Dantzig Selector satisfy sparsity inequalities under weaker assumptions than those required for the "original" LASSO.

\end{abstract}


\section{Introduction}

In many modern applications, a statistician often have to deal with very large datasets. They may involve a large number $p$ of covariates, possibly larger than the sample size $n$. Let an observation be a pair instance-label.
In this paper, we tackle such high dimensional settings which moreover involve a large amount of unlabeled data (say $m$ instances) in addition to the $n$ labeled observations.

In contrast to inductive or supervised methods, transductive procedures exploit the knowledge of the unlabeled data to improve prediction.
It is argued in the semi-supervised learning literature (see for example \cite{semi-sup} for a recent survey) that taking into account the information
on the design given by the new additional instances has a stabilizing effect on
the estimator. In transductive methods, we furthermore take into account the objective
of the statistician: estimation of the value of the regression function only
on the set of unlabeled data; see Vapnik \cite{Vapnik} for a pioneer work.
To leverage unlabeled data, the transductive or semi-supervised methods exploit the geometry of the marginal distribution.
Initially introduced in the classification framework, the good performance of these methods has been observed in several practical fields.
From a theoretical point of view we refer to the transductive version of SVM
\cite{Vapnik98,Joach99TSVM,Chap05Trans,WangTransSVM07}, to the study of classifiers
under the clustering assumption \cite{Blum98Trans,Nigam98ENTrans,Collin99,Zhou99TransGraph,Zhu99Trans,Ando05Trans}, and
to transductive versions of the Gibbs estimators \cite{manuscrit} among many others.
According to the applications arrays, see for instance
the detection of spam email \cite{Blum98Trans,Amini03Trans,Balcan05Trans}, genetics applications \cite{XiaoTrans05}, but also the well-known Netflix challenge.

In this paper we focus on the linear regression model in the case $p>n$.
In this setting dimension reduction is a major issue and can be performed through the selection of a small amount of relevant covariates.
Numerous inductive or supervised methods have been proposed in the literature, ranging from the classical information criteria such as $\mathop{\rm AIC}$ \cite{aic} and $\mathop{\rm BIC}$ \cite{bic} to the more recent $ \ell_{1}$-regularized methods, as the LASSO \cite{Tibshirani-LASSO}, the Dantzig Selector \cite{Dantzig},
the non-negative garrote \cite{garrotte}. We also refer to \cite{KoltchDant,Koltchl1plus,VanGpLass,VandeGeerSparseLasso,ArnakTsyb,ChriMo7GpLass}
for related works.
Such regularized regression methods have recently witnessed several developments due to the attractive feature of computational feasibility, even for high dimensional data ({\it i.e.}, when the number of covariates $p$ is large).
Formally we assume
\begin{equation}
\label{eq_depart}
	y_{i}
	=
	x_{i} \beta^* + \varepsilon_{i},
	\quad \quad i=1,\ldots,n,
\end{equation}
where the design $x_{i}=(x_{i,1},\ldots,x_{i,p}) \in \mathbb{R}^p$ is deterministic,
$\beta^*=(\beta^*_1,\ldots,\beta^*_p)' \in \mathbb{R}^p$ is the unknown parameter
and $\varepsilon_1,\ldots,\varepsilon_n$ are i.i.d. centered Gaussian random
variables with known variance $\sigma^{2}$.
Let $X$ denote the matrix with $i$-th line equal to $x_i$, and let $X_{j}$ denote its $j$-th column, with $i\in\{1,\ldots,n\}$ and $j\in\{1,\ldots,p\}$.
In this way, we can write
$$
X
=
(x'_1 , \ldots, x'_n )'
= (X_1,\ldots,X_p).
$$
For the sake of simplicity, we will assume that the observations are normalized in such a way that $X_{j}'X_{j}/n=1$. We denote by $Y$ the vector $Y=(y_{1},\ldots,y_{n})'$.

Let $x_{n+1},\ldots x_{m}$ be observed unlabeled instances with $x_{i}\in\mathbb{R}^p$ for $n+1\leq i \leq m$ (with $m>n$).
Let moreover $Z=(x'_{1} , \ldots, x'_{m} )'$.

For all $\alpha \geq 1$ and any vector $v\in\mathds{R}^{d}$, we set $ \|\cdot \|_{\alpha}$, the $\ell_{\alpha}$-norm: $ \|v\|_{\alpha} = (|v_{1}|^{\alpha}+\ldots+|v_{d}|^{\alpha})^{1/\alpha}$.
In particular $\|\cdot\|_{2}$ is the euclidean norm.
Moreover for all $d$-dimensional vector $v$ with $d\in\mathds{N}$, we use the notation $ \|v\|_{0} = \sum_{i=1}^{d} \mathds{1}(v_{i}\neq 0). $

From a theoretical point of view, Sparsity Inequalities (SI) have been proved for the regularized estimators mentioned above under different assumptions, in the inductive setting only,
{\it i.e.}, without the knowledge of the matrix $Z$.
That is upper bounds of order of $\mathcal{O}\left(\sigma^{2} \|\beta^*\|_{0} \log(p)/n\right)$ for the errors $(1/n)\|X\hat{\beta} - X\beta^{*}\|_{2}^{2}$
and $\|\hat{\beta} - \beta^{*}\|_{2}^{2}$ have been derived, where $\hat{\beta}$ is one of those estimators.
Such bounds involve the number of non-zero coordinates in $\beta^{*}$ (multiplied by $\log(p)$), instead of dimension $p$.
Such bounds guarantee that under some assumptions, $X \hat{\beta}$ and $\hat{\beta}$ are good estimators of $X\beta^{*}$ and $\beta^{*}$ respectively.
For the LASSO, these SI are given for example in \cite{BTWAggSOI,Lasso3}, whereas \cite{Dantzig,Lasso3} provided SI for the Dantzig Selector.
On the other hand, Bunea \cite{Bunea_consist} established conditions which ensure that the LASSO estimator and $\beta^{*}$ have the same null coordinates.
Analog results for the Dantzig Selector can be found in~\cite{KarimNormSup}.
An important issue when we establish these theoretical results is the assumption  that is needed on the Gram matrix $n^{-1}X'X$.
We refer to \cite{VandeGeerConditionLasso09} for a nice overview of these assumptions.

In this paper we are interested in the estimation of $Z\beta^{*}$: namely, we care about predicting what would be the labels attached to the additional $x_{i}$'s.
Hence we develop transductive versions of the LASSO or the Dantzig Selector to tackle the problem of estimating the vector $Z\beta^*$ in the high dimensional setting.
According to \cite{Vapnik} this estimator should differ from an estimator tailored for the estimation of $\beta^{*}$ or $X\beta^{*}$ like the LASSO.
Indeed, a naive plug-in method would be to build an estimator $\hat{\beta}(X,Y)$ and then to compute $Z\hat{\beta}(X,Y)$ to estimate $Z\beta^*$.
We rather consider here an approach where the estimators $\hat{\beta}(X,Y,Z)$ exploit the knowledge of $Z$, and finally compute $Z\hat{\beta}(X,Y,Z)$.
These transductive procedures observe several interesting properties:
\begin{itemize}	
	\item they take advantage of the unlabeled points to satisfy SIs with weaker assumptions on the Gram matrix than those required for the usual inductive methods (cf. the examples of Section~\ref{exam});
	\item they perform well in practice compared to the inductive methods in most of the situations.
	This is illustrated by a comparison between the performance of the Transductive LASSO and the LASSO.
\end{itemize}

Let us mention that the study established in this paper consists in a generalization of the LASSO and the Dantzig Selector. They actually do not only consider the transductive setting since they can be adapted to other objectives desired by the statistician.
Indeed, the estimators depend on a $q\times p$ matrix $A$, with $q\in\mathds{N}$, whose choice allows to consider the problem of the estimation of $A\beta^{*}$.
In this way, the estimation of $Z\beta^{*}$ appears as a particular case.

The rest of paper is organized as follows.
In the next section, we give the definition of the considered estimators.
We then display, in Section~\ref{simu}, a set of experiments that show how the Transductive estimators can improve on the LASSO and the Dantzig Selector in many applications.
A non-asymptotic study is provided in Section~\ref{thms} whereas all the proofs of the theorems are postponed to Section~\ref{proofs}.

\section{Definition of the estimators}

Here we define the family of estimators we consider in the sequel and more specifically the Transductive LASSO and the Transductive Dantzig Selector.

\subsection{Definitions}
\label{sec:Def}

Let us first remind that the LASSO estimate can be defined by
\begin{equation}
\label{LASSO-def1}
	\hat{\beta}^{L}_{\lambda} 
	=
	 \arg\min_{\beta \in \mathds{R}^{p}} \left\{\left\|Y-X\beta\right\|_{2}^{2} + 2\lambda \|\beta\|_{1} \right\},
\end{equation}
where $\lambda$ is a positive tuning parameter. Let us make simple remarks to optimize comprehensibility of the paper and the notation inside. Since $Y=X\beta^{*}+\varepsilon$, the response vector $Y$ can be seen as an estimator for $X\beta^{*}$.
Then let us write $Y=\widehat{X\beta^{*}}$ to convey this fact.
Actually, if $\sigma\simeq 0$, $Y$ could even be a good estimator.
However, in the general case, it is not expected to be a particularly interesting estimator and then $Y=\widehat{X\beta^{*}}$ is only used as a preliminary estimator of the vector $X\beta^{*}$.
Based on this preliminary estimate, the LASSO defined by~\eqref{LASSO-def1} ensures in the case where $\beta^{*}$ is sparse, a good estimation of $X\beta^{*}$ by $X\hat{\beta}^{L}_{\lambda}$ (cf. \cite{Lasso3} for instance).

Let us now generalize the previous comments. Let $A\beta^{*}$ be a quantity of interested for a general (and given) $q\times p$ matrix $A$ with $q\in\mathbb{N}^*$. Then an analog of the LASSO estimator~\eqref{LASSO-def1} can be given by the following definition.
\begin{dfn}[The Transductive LASSO]
\label{dfntlasso}
	Let $\widehat{A\beta^{*}}$ be a preliminary estimator (that can be a very poor estimator) of $A\beta^{*}$ and define
	$$
	\hat{\beta}_{A,\lambda}
	=
	\arg\min_{\beta \in \mathds{R}^{p}} \left\{\left\|\widehat{A\beta^{*}}-A\beta\right\|_{2}^{2} + 2\lambda \|\beta\|_{1} \right\}.
	$$
	In particular, when $A=\sqrt{n/m}Z$, the estimator $\hat{\beta}_{A,\lambda} $
            is the Transductive LASSO.
\end{dfn}
The Dantzig Selector is defined as:
\begin{equation}
\label{DS-def}
	\tilde{\beta}^{DS}_{\lambda}
	=
	\left\{
	\begin{array}{l}
		\arg\min_{\beta \in \mathds{R}^{p}} \left\|\beta\right\|_{1}
		\\
		\\
		s. t. \left\|X'(Y-X\beta)\right\|_{\infty} \leq \lambda,
	\end{array}
	\right.
\end{equation}
where $\lambda$ is a positive tuning parameter. In the same way as for the LASSO estimator, we underline the role of $Y=\widehat{X\beta^{*}}$ and propose the following definition.

\begin{dfn}[The Transductive Dantzig Selector]
\label{dfntds}
	Let $\widehat{A\beta^{*}}$ be a preliminary estimator of $A\beta^{*}$ and define
	$$
	\tilde{\beta}_{A,\lambda}
	 =
	\left\{
	\begin{array}{l}
		\arg\min_{\beta \in \mathds{R}^{p}} \left\|\beta\right\|_{1}
		\\
		\\
		s. t. \left\|A'(\widehat{A\beta^{*}}-A\beta)\right\|_{\infty} \leq \lambda.
	\end{array}
	\right.
	$$
	In particular, when $A=\sqrt{n/m}Z$, the estimator $\tilde{\beta}_{A,\lambda} $ is the Transductive Dantzig Selector.
\end{dfn}

In Definitions~\ref{dfntlasso} and~\ref{dfntds}, the matrix $A$ is not specified. Hence, we cover here a general objective. However, we mainly focus in this paper on the transductive setting. Then our principal study deals with the estimation of $Z\beta^*$. In other works, this means that we set in the above definitions $A=\sqrt{n/m}Z$, the unlabeled data matrix
(with a normalization term $\sqrt{m/n}$ that is here for the sake of convenience, see Section~\ref{thms}).
An important issue in this paper is also the preliminary estimator $\widehat{A\beta^{*}}$ that should be used.
An explicit condition on this estimator can be found in Section~\ref{thms}.
It ensures the good theoretical performance for $\hat{\beta}_{A,\lambda}$ and
$\tilde{\beta}_{A,\lambda}$.
Let us now propose some examples of preliminary estimators.

\begin{exms}
\label{ex:PreliminaryEstimator}
	i) The most simple idea is to estimate $\beta^{*}$ by the least square estimator.
	Even if $p>n$, that implies that $(X'X)$ is not invertible, we can choose any pseudo-inverse $\widetilde{(X'X)}^{-1}$ of $(X'X)$ and use in this way
	$$
	\widehat{A\beta^{*}}
	=
	A \widetilde{(X'X)}^{-1} X' Y ,
	$$
	as preliminary estimator of $A\beta^{*}$.
	Remark that if $\Ker(A)\subset\Ker(X)$, this quantity is uniquely defined (it does not depend on the choice of the particular pseudo-inverse $\widetilde{(X'X)}^{-1}$).
	We will see later that in this case, we may have theoretical guarantees for the performance of $\hat{\beta}_{A,\lambda}$ and $\tilde{\beta}_{A,\lambda}$.\\
	\noindent ii) One may think about more sophisticated regularization procedures, based for instance on 
	$$
	\widehat{A\beta^{*}}
	=
	A (\gamma A'A + X'X)^{-1} X' Y,
	$$
	for (a small) $\gamma>0$ when the matrix $A'A$ is invertible (in the idea of ridge regression).\\
	\noindent iii) Finally, practitioners may prefer to use as a preliminary estimator something known to work well in practice, like:
	$$
	\widehat{A\beta^{*}}
	=
	A \hat{\beta}^{L}_{\lambda'} ,
	$$
	with $\lambda'\geq 0$. We pay a particular attention to this initial estimator in the rest of the paper.
\end{exms}

\subsection{A discussion on the matrix $A$}
\label{Sec:DiscussA}

These novel estimators depend on two tuning parameters.
First $\lambda>0$ is a regularization parameter, it plays the same role as the tuning parameter involved in the LASSO and will be discussed in our simulations and our theoretical results.
The second one is the matrix $A$, that allows to adapt the estimator to the objective of the statistician.
In this paper, we are mainly interested in the following objective:
\begin{itemize}
	\item {\bf transductive objective:} the estimation of $Z\beta^{*}$, by $Z\hat{\beta}_{A,\lambda}$ or $Z\tilde{\beta}_{A,\lambda}$ with $A=\sqrt{n/m}Z$. Note that in the case
$n<p<m$, it is possible that the matrix $Z'Z$ is invertible, while $X'X$ may not.
\end{itemize}
Other choices of the matrix $A$ are possible and help to deal with other objectives. We display here two additional feasible and well-known choices:
\begin{itemize}
	\item {\bf denoising objective:} the estimation of $X\beta^{*}$, that is a denoised version of $Y$.
	For this purpose, we consider the estimator $\hat{\beta}_{A,\lambda}$, with $A=X$.
	In this case, if we keep $Y$ as our preliminary estimator of $X\beta^{*}$, our estimators are exactly the LASSO and the Dantzig Selector, so this case is not of particular interest in this paper;
	\item {\bf estimation objective:} the estimation of $\beta^{*}$ itself, by $\hat{\beta}_{A,\lambda}$, with $A=\sqrt{n}\mathbf{I}_{p}$ where $\mathbf{I}_{p}$ is the identity matrix of size $p$.
\end{itemize}
Thanks to the unifying notation $A$, the theoretical performance of the estimators based on the above choices are considered in the same time in Section~\ref{thms}.
However, we mention that the main contribution relates to the transductive objective.

\section{Experimental results}
\label{simu}

In this section we compare the empirical performance of the Transductive LASSO and the LASSO estimators on simulated and real datasets according to the transductive task.
We consider both low and high dimensional simulated data.
The real dataset comes from a genetic study, devoted to learn the complex combinatorial code underlying gene expression.
More precisions are given in Section~\ref{sec:RealData}.
In this dataset, there are $p = 666$ predictor variables and the total number of available labeled data is $2587$.
The conclusions of our experiments is that the transductive LASSO outperforms the LASSO estimator in most settings and specifically when the variables are correlated.

\subsection{Implementation}

Since the paper of Tibshirani \cite{Tibshirani-LASSO}, several effective algorithms to compute the LASSO have been proposed and studied (for instance Interior Points methods \cite{interior}, LARS \cite{Efron-LARS}, Pathwise Coordinate Optimization \cite{PCO}, Relaxed Greedy Algorithms \cite{Barron2}).
For the Dantzig Selector, a linear method was proposed in the first paper \cite{Dantzig}.
The LARS algorithm was also successfully extended in \cite{DASSO} to compute
the Dantzig Selector.\\
\noindent Note that these methods allow to compute our estimators $\hat{\beta}_{A,\lambda}$ and $\tilde{\beta}_{A,\lambda}$ as they just appear as the LASSO and the Dantzig Selector computed on modified data.
Namely, after the computation of the preliminary estimator $\widehat{A\beta^{*}}$ of $A\beta^{*}$, the transductive estimator $\hat{\beta}_{A,\lambda}$ is obtained as a usual LASSO solution where the usual data $(X,Y)$ are replaced by $(A,\widehat{A\beta^{*}})$.

We use the version of the Transductive LASSO proposed in Section~\ref{dfntlasso} based on the LASSO as preliminary estimator.
That is, we set $A=Z$ in Definition~\ref{dfntlasso} and refer to Example~\ref{ex:PreliminaryEstimator}-iii) for the definition of the preliminary estimator $\widehat{A\beta^{*}}$.
In other words, for a given $\lambda_{1}$, we first compute the LASSO estimator $\hat{\beta}_{X,\lambda_{1}} = \hat{\beta}^{L}_{\lambda_1}$.
In this way we have $\widehat{Z\beta^{*}} = Z\hat{\beta}^{L}_{\lambda_1}$.
However, recalling that $Z=(x_1',\ldots , x_{n}',\ldots,x_m')'$, it is worth noting that we should keep the $n$ first components of $\widehat{Z\beta^{*}}$ equal to $Y$.
Indeed the $n$ first components correspond to the labeled samples and then do not require to be replaced by an estimation. 
Given this adjustment, the Transductive LASSO is given by
\begin{equation*}
\hat{\beta}^{TL}(\lambda_{1},\lambda_{2})= \left\{
\begin{array}{l}
\arg\min_{\beta \in \mathds{R}^{p}} \frac{n}{m} \left\|Z\beta\right\|_{2}^{2}
\\
\\
s. t. \left\|\frac{n}{m} Z'(Z \hat{\beta}^{L}_{\lambda_1} - Z\beta) \right\|_{\infty} \leq \lambda_{2},
\end{array}
\right.
\end{equation*}
for a given $\lambda_{2}$ (cf. Section~\ref{thms} for a theoretical study of this estimator).
Let us mention that the good performance of this estimator are stated in Theorem~\ref{thmpreliminarylasso} under some assumptions.
We compare this two steps procedure with the procedure obtained using the usual only LASSO $\hat{\beta}_{\lambda}^{L}=\hat{\beta}_{X,\lambda}$ for a given $\lambda$ that may differ from $\lambda_{1}$.
In both cases, the solutions are computed using the {\it glmnet}\footnote{The algorithm is implemented with R and can be found in the web page: http://cran.r-project.org/web/packages/glmnet/index.html} package, introduced by Friedman et al., to provide the LASSO solution, between others.
We compute $\hat{\beta}_{\lambda}^{L}$ and $\hat{\beta}^{TL}(\lambda_{1},\lambda_{2})$ for $(\lambda,\lambda_{1},\lambda_{2})\in\Lambda^3$ where $\Lambda$ is some grid defined in a data driven way by the glmnet algorithm.
In all the experiments, we choose the best tuning parameters.
That is, the choice is based on the truth. In other words, we only compare
the oracle in our family of estimators.
This way to select the tuning parameter is
even suitable in our real data experiments. Indeed, the initial data we get consist only in labeled data.
We then hide many responses values and construct the estimators without their knowledge.
Finally the best estimators (tuning parameters) are chosen based on those hidden responses.

\subsection{Synthetic data}

The comparison between the LASSO and the Transductive LASSO is made through the study of the distribution of
$$ PERF(Z) = \frac{\min_{(\lambda_{1},\lambda_{2})\in\Lambda^{2}} \| Z(\hat{\beta}^{TL}(\lambda_{1},\lambda_{2})-\beta^*)
      \|_{2}^{2} }{ \min_{\lambda\in\Lambda} \| Z(\hat{\beta}_{\lambda}^{L}-\beta^*)
               \|_{2}^{2}}, $$
over $100$ replications for each experiment.
Since the LASSO is a special case of the Transductive LASSO (with $\lambda_2=0$), $ PERF(Z)$ belong to $[0,1]$.
This ratio measures the improvement made by the Transductive LASSO according to the transductive objective.
The smaller $ PERF(Z)$ is, the more attractive the use of the Transductive LASSO is.
Analog study can be considered to compare the performance of the Transductive LASSO and the LASSO in term of the denoising and the estimation tasks respectively based on the ratio
$$ PERF(X) = \frac{\min_{(\lambda_{1},\lambda_{2})\in\Lambda^{2}} \| X(\hat{\beta}^{TL}(\lambda_{1},\lambda_{2})-\beta^*)
      \|_{2}^{2} }{ \min_{\lambda\in\Lambda} \| X(\hat{\beta}_{\lambda}^{L} -\beta^*)
               \|_{2}^{2}}, $$
and
$$ PERF(I) = \frac{\min_{(\lambda_{1},\lambda_{2})\in\Lambda^{2}} \| \hat{\beta}^{TL}(\lambda_{1},\lambda_{2})-\beta^*
      \|_{2}^{2} }{ \min_{\lambda\in\Lambda} \| \hat{\beta}_{\lambda}^{L}-\beta^*
               \|_{2}^{2}}.$$\\

\noindent{\bf Data description.} We consider several simulations from the linear regression model
$$
y_{i}
=
x_i \beta^* + \varepsilon_i,
$$
for $i\in\{1,\ldots,n\}$, $\beta^* \in\mathds{R}^{p}$ and the $\varepsilon_{i}$ are i.i.d. $\mathcal{N}(0,\sigma^2)$.
The design matrix comes from a centered multivariate normal distribution with covariance structure $Cov(X_j ;X_k) = \rho^{-|j-k|} $ with $\rho\in]0,1[$.
Dimension $p$, sample sizes $(n,m)$, noise level $\sigma$ and correlation parameter $\rho$ are left free.
They will be specified during the experiments in order to check the robustness of the results.
The regression vector $\beta^*$ is $s$-sparse where $s \geq 1$ is an integer and corresponds to the sparsity index.
That is, $\beta^*$ consists in $s$ non-zero components.
Since the LASSO and the Transductive LASSO do not take care of the ordering of the variables let us define $\beta^*$ such as its $s$ first components are non-zero and equal $5$.\\
\noindent Our study covers several combinations of the parameters $p$, $s$, $(n,m)$, $\rho$ and $\sigma^2$.
In the next paragraph, we examine the performance of each estimator according to the value of the regularization parameters.\\

\noindent {\bf Results.} We consider separately the low and the high dimensional case.

\vspace{0.5cm}

\begin{table}[t]
\caption{Evaluation of the mean ($Mean$), the median ($Med$) and the quantile $Q_3$ of order $0.3$ of the quantities $PERF(Z)$, $PERF(X)$ and $PERF(I)$, when the methods are used in the artificial low dimensional case $p < n$.}
\label{tab:frqVSsigma}
\begin{center}
\begin{footnotesize}
\begin{sc}
\begin{tabular}{|l|c|c|c|c||c|c|c||c|c|c||c|c|r|}
\cline{6-14}
 \multicolumn{5}{c|}{} & \multicolumn{3}{c||}{$PERF(Z)$}    & \multicolumn{3}{c||}{$PERF(X)$}    &\multicolumn{3}{c|}{$PERF(I)$}    \\
\hline
$p$   & $s$   &  $(n,m)$ & $\rho$ & $\sigma^2$ &  $Mean$  & $Med$ &  $Q_3$  & $Mean$    & $Med$ &  $Q_3$  & $Mean$  & $Med$ &  $Q_3$  \\
\hline
\hline
 $8$ & $1$ & $(10,30)$ &  $0.1$  & $1$ &  $0.64$ &   $0.77$  & $0.47$  & $0.62$ & $0.70$& $0.44$ & $0.64$  & $0.75$ & $0.42$  \\
\hline
 $8$ & $1$ & $(10,30)$ &  $0.9$  & $25$ &  $0.57$ & $0.54$ & $0.11$  & $0.63$  & $0.47$ & $0.11$    & $0.66$ & $0.61$  &  $0.10$ \\
\hline
$ 50$ & $1$ & $(60,80)$&  $0.1$   & $1$  &  $0.73$  & $0.81 $ &  $0.59$ & $0.74$  &  $0.82$ & $0.62$  & $0.72$ & $0.77$ &  $0.60$    \\
\hline
$50$ & $1$ & $(60,200)$ & $0.9$   & $1$ &  $0.70$  & $0.86$ & $0.63$ &  $0.69$ & $0.83$ &  $0.60$ & $0.68$  & $0.80$ & $0.56$       \\
\hline
\hline
$8$& $3$ & $(10,30)$ &  $0.1$  & $1$ & $0.85$   & $0.91$ & $0.84$  &  $0.83$ & $0.89 $ & $0.81$    & $0.83$ & $0.93$ &  $0.78$    \\
\hline
$8$& $3$ & $(10,30)$ &  $0.9$  & $100$ & $0.74$   & $0.84$ & $0.70$  &  $0.71$ & $0.80 $ & $0.59$    & $0.75$ & $0.79$ &  $0.61$    \\
\hline
$8$& $3$ & $(10,100)$ &  $0.5$  & $1$ & $0.90$   & $0.98$ & $0.91$  &  $0.89$ & $0.99 $ & $0.88$    & $0.87$ & $0.95$ &  $0.85$    \\
\hline
$8$& $3$ & $(10,100)$ &  $0.9$  & $25$ & $0.75$   & $0.88$ & $0.64$  &  $0.72$ & $0.80 $ & $0.60$    & $0.74$ & $0.82$ &  $0.68$    \\
\hline
$50$& $20$ & $(100,120)$ & $0.1$   & $1$ &  $0.98$  &  $1$ & $0.98$ &  $0.98$ & $1$  & $0.98$ & $0.98$ & $1$  & $0.98$ \\
\hline
$50$ & $20$ & $(100,120)$ & $0.9$   & $1$ &  $0.76$  & $0.75$ & $0.68$ &  $0.58$ & $0.56$ &  $0.51$ & $0.96$  &  $1$& $1$       \\
\hline
\end{tabular}
\end{sc}
\end{footnotesize}
\end{center}
\vskip -0.1in
\end{table}

\noindent{\it The low dimensional case.} Several examples have been studied and we illustrate the performance of the methods in this case through some specific experiments.
The main parameter seems to be the sparsity index $s$.
More precisely, the behavior of the Transductive LASSO compared to the LASSO is related to how large the sparsity index is in comparison to the dimension $p$.
This is illustrated in Table~\ref{tab:frqVSsigma}, where the two cases are separated by two horizontal lines.
Hence, when $s$ is small, we notice a good improvement while using the Transductive LASSO instead of the LASSO estimator.
This is displayed in lines~1 to~3 of Table~\ref{tab:frqVSsigma}, where all of the ratios $PERF(Z)$, $PERF(X)$ and $PERF(I)$ are most of the time between $0.5$ and $0.8$.
We remark also that the performance of the Transductive LASSO are even better when the parameters $\rho$ and $\sigma$ increase (line~2). 
On the other hand, when the sparsity index $s$ is large (with respect to the dimension $p$), it turns out that the Transductive LASSO does not improve enough the LASSO estimator (lines~7 and~9) when $p$ is large, whereas it is still satisfying for small values of $p$.
By poor improvement, we mean that $\min_{(\lambda_{1},\lambda_{2})\in\Lambda^{2}}
\| Z(\hat{\beta}^{TL}(\lambda_{1},\lambda_{2})-\beta^*) \|_{2}^{2}$ is not far from $\min_{(\lambda_{1},0)\in\Lambda^{2}} \| Z(\hat{\beta}^{TL}(\lambda_{1},0)-\beta^*) \|_{2}^{2}$ which coincides with the LASSO error $
\min_{\lambda\in\Lambda} \| Z(\hat{\beta}_{\lambda}^{L} -\beta^*)\|_{2}^{2}$.
An important observation is that even when $s$ is large and in the case where the variables are highly correlated, that is when $\rho$ is large, the Transductive LASSO can be a good alternative to the LASSO estimator (lines~6, 8 and~10).
This observation is true for both of the transductive error ratio $PERF(Z)$ and the denoising error ratio $PERF(X)$.
On the other hand, even with high correlations between variables, the Transductive LASSO does not make the estimation error ratio $PERF(I)$ better in this last situation.\\
In the low dimensional case, it seems that increasing the number $m$ of unlabeled data does not lead to an improvement of the performance of the Transductive LASSO.
This is for instance displayed in line~4 of Table~\ref{tab:frqVSsigma} where  $m=200$ or in lines~7 and~8 of the same table, where $m=100$.
Finally, note that when $n$ is large, the LASSO estimator behaves in a good way and it becomes difficult to improve it thanks to the Transductive LASSO.

\vspace{0.5cm}

\begin{table}[t]
\caption{Evaluation of the mean ($Mean$), the median ($Med$) and the quantile $Q_3$ of order $0.3$ of the quantities $PERF(Z)$, $PERF(X)$ and $PERF(I)$, when the methods are used in the artificial high dimensional case $p \geq n$.}
\label{tab:ResultGdDim}
\begin{center}
\begin{scriptsize}
\begin{sc}
\begin{tabular}{|l|c|c|c|c||c|c|c||c|c|c||c|c|r|}
\cline{6-14}
 \multicolumn{5}{c|}{} & \multicolumn{3}{c||}{$PERF(Z)$}    & \multicolumn{3}{c||}{$PERF(X)$}    &\multicolumn{3}{c|}{$PERF(I)$}    \\
\hline
$p$ & $s$    &  $(n,m)$ & $\rho$ & $\sigma^2$ &  $Mean$  &  $Q_3$ & $Med$ & $Mean$ & $Med$  &  $Q_3$  & $Mean$  & $Med$ &  $Q_3$  \\
\hline
\hline
$10$ & $8$ & $(5,15)$ &  $0.1$  & $1$ &  $0.73$  & $0.77$  & $0.64$ & $0.33$ & $0.23$ & $0.14$ & $0.73$  & $0.74$ & $0.65$  \\
\hline
$10$& $8$  & $(5,15)$ &  $0.9$  & $1$ &  $0.78$  &   $0.84$  & $0.71$ & $0.45$ & $0.42$ & $0.29$ & $0.78$  & $0.80$ & $0.73$  \\
\hline
$10$ & $8$ & $(5,15)$ &  $0.9$  & $100$ &  $0.79$  & $0.82$  & $0.71$  & $0.70$& $0.77$ & $0.55$    & $0.79$  & $0.81$&  $0.72$ \\
\hline
$1000$ & $50$ & $(20,60)$ &  $0.1$  & $1$ &  $0.94$  & $1$  & $1$  & $0.95$& $1$ & $1$    & $0.97$  & $1$&  $1$ \\
\hline
$1000$ & $50$ & $(20,60)$ &  $0.9$  & $1$ &  $0.60$  & $0.53$  & $0.36$  & $0.44$& $0.48$ & $0.33$    & $0.78$  & $0.82$&  $0.65$ \\
\hline
\hline
$ 1000$& $50$ & $(100,200)$ &  $0.1$  & $1$ & $0.98$   & $1$ & $1$  &  $0.87$ & $0.89$ & $0.79$    & $0.99$  & $1$&  $1$    \\
\hline
$1000$ & $50$ & $(100,200)$&  $0.9$   & $1$  &  $0.49$  & $0.45$ & $0.38$ & $0.35$  & $0.32$ &   $0.25$ & $0.90$  & $0.97$ &  $0.86$    \\
\hline
$1000$& $50$ & $(100,200)$ & $0.9$   & $25$ &  $0.50$  & $0.46$  &  $0.39$ &  $0.43$ & $0.40$ & $0.30$ & $0.86$  & $0.91$ &  $0.80$       \\
\hline
\end{tabular}
\end{sc}
\end{scriptsize}
\end{center}
\vskip -0.1in
\end{table}

\begin{figure}[t]
\vskip -0.2in
\includegraphics[width=2in,height=1.5in]{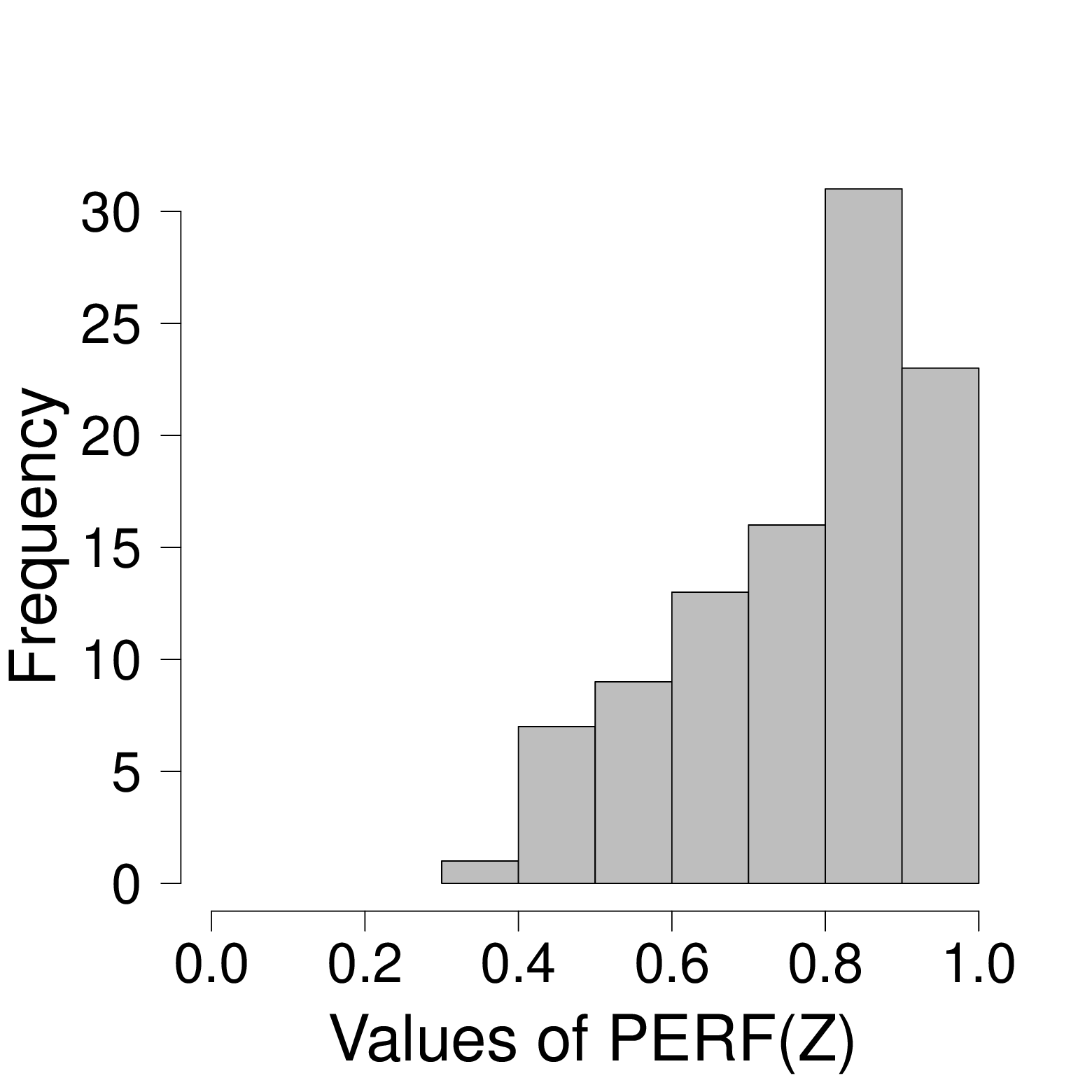}
\includegraphics[width=2in,height=1.5in]{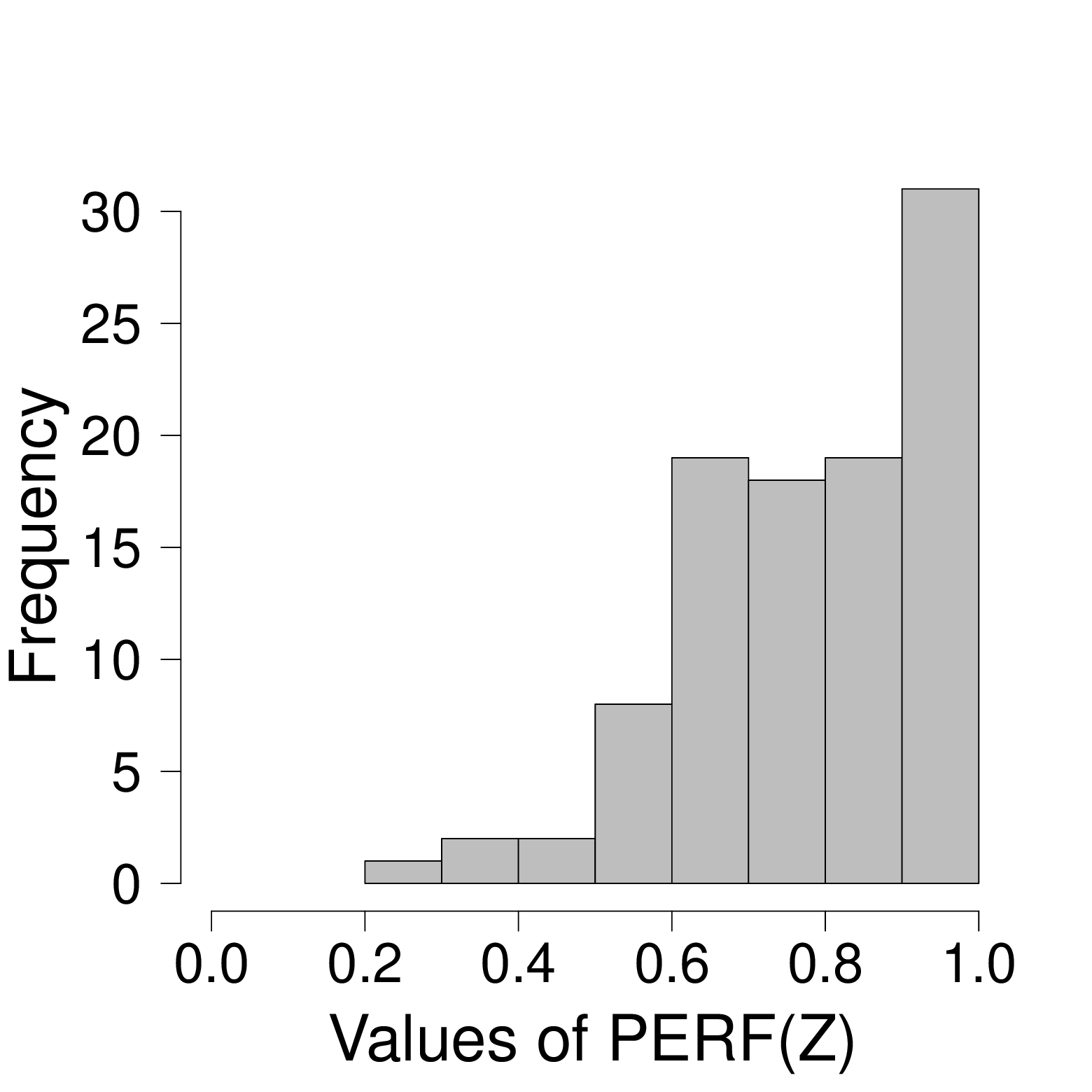}
\includegraphics[width=2in,height=1.5in]{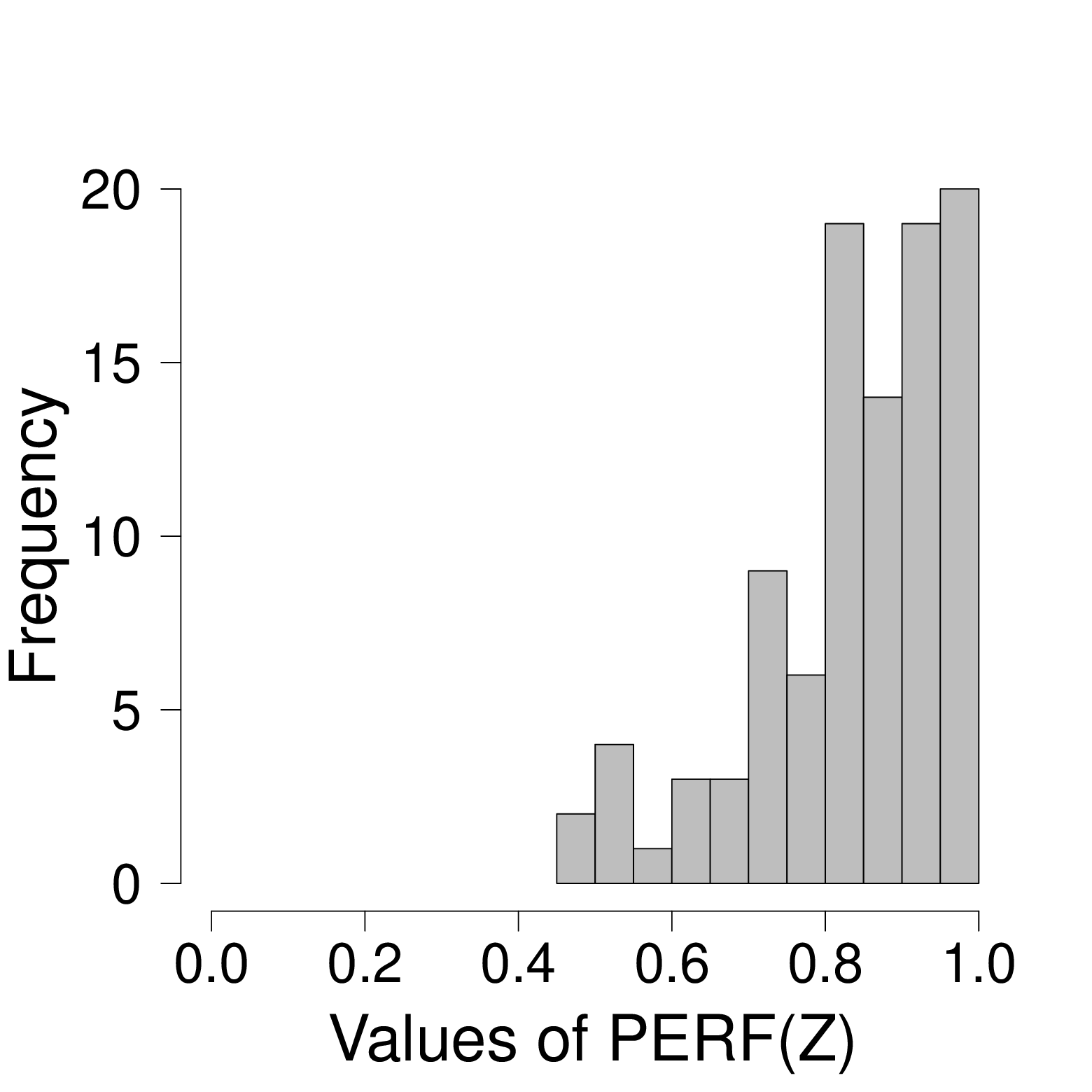}
\\
\includegraphics[width=2in,height=1.5in] {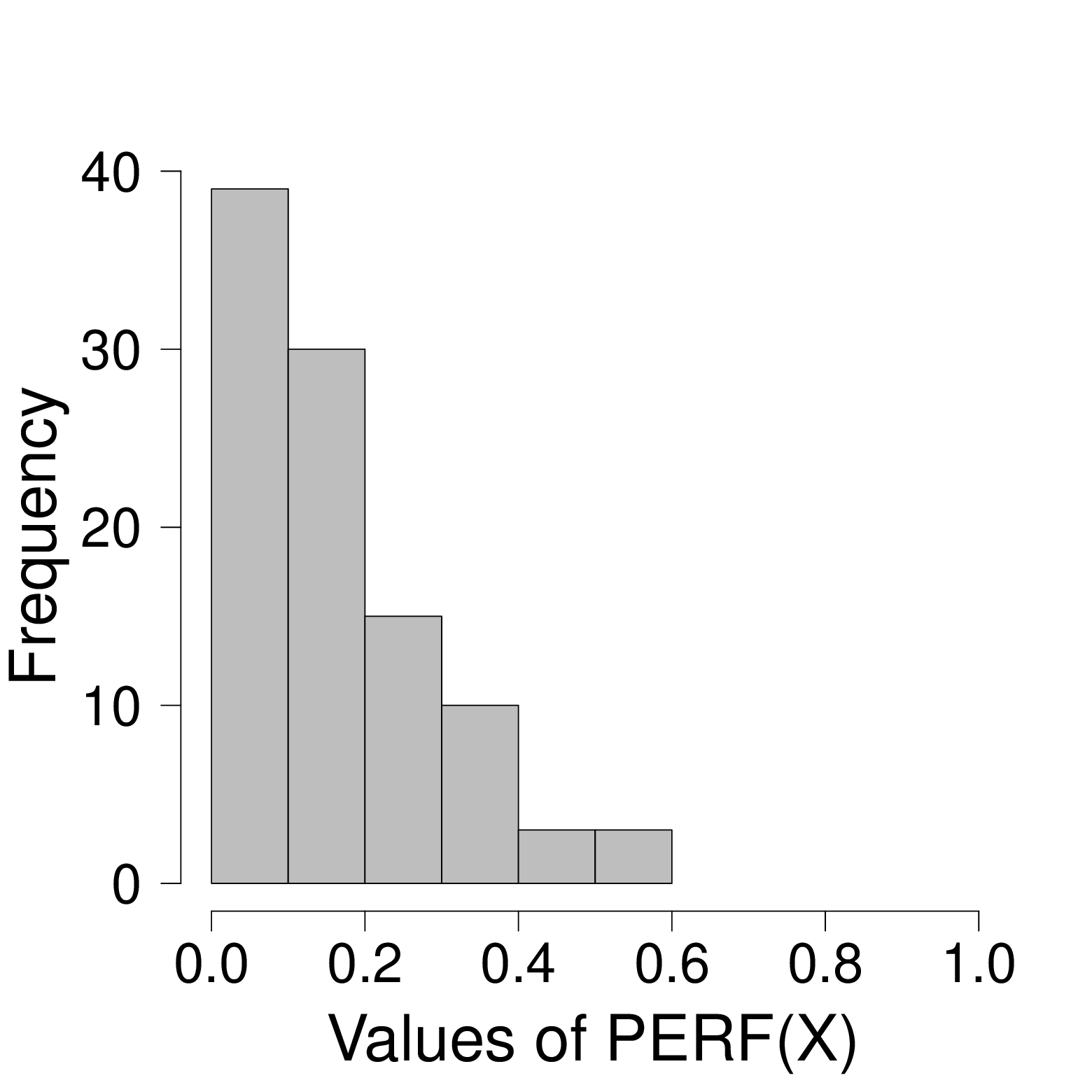}
\includegraphics[width=2in,height=1.5in] {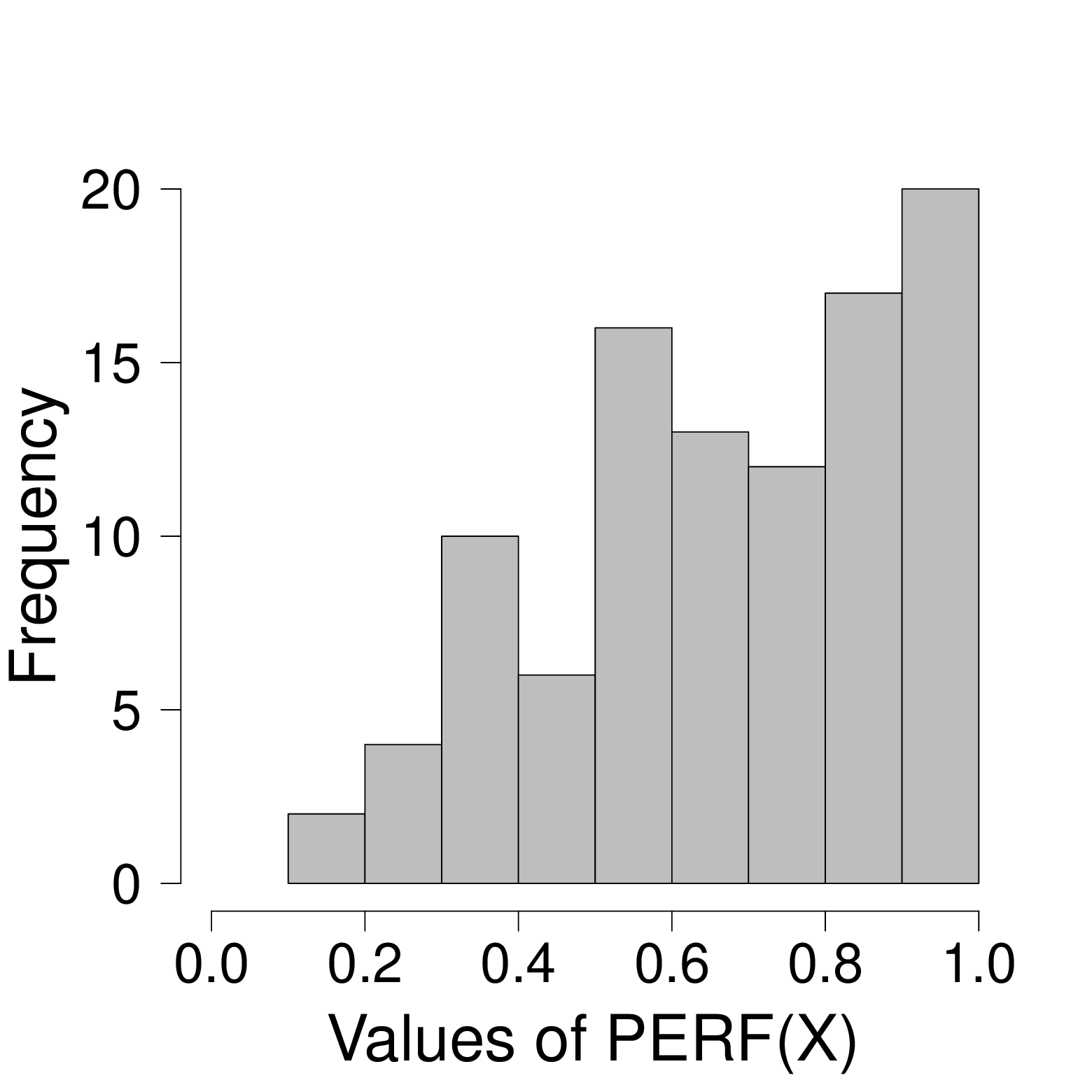}
\includegraphics[width=2in,height=1.5in] {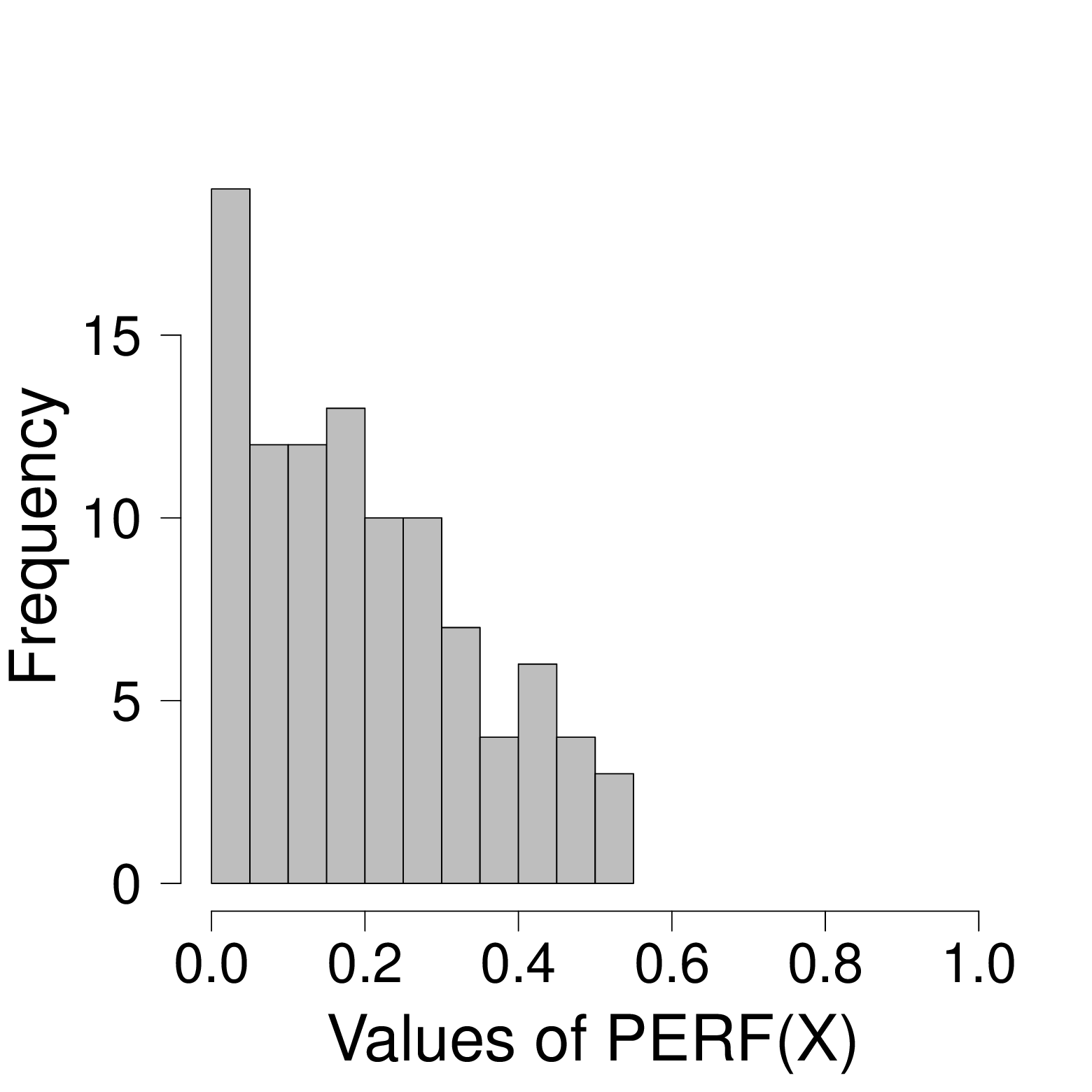}
\vglue-40pt
~
\begin{center}
\begin{minipage}[t]{0.90\textwidth}
\caption{\footnotesize Study of the distribution of $PERF$ when $p=50$, $s=20$ and $(n,m)  = (10,20)$. Top: study of $PERF(Z)$; Bottom: study of $PERF(X)$. From left to right: $\rho = 0.1$ and $\sigma^2 = 1$; $\rho = 0.1$ and $\sigma^2 = 100$; $\rho = 0.9$ and $\sigma^2 = 1$.}
\label{fig:histMultiplp50s20}
\end{minipage}
\end{center}
\end{figure}

\begin{figure}[t]
\vskip -0.2in
\includegraphics[width=1.5in]{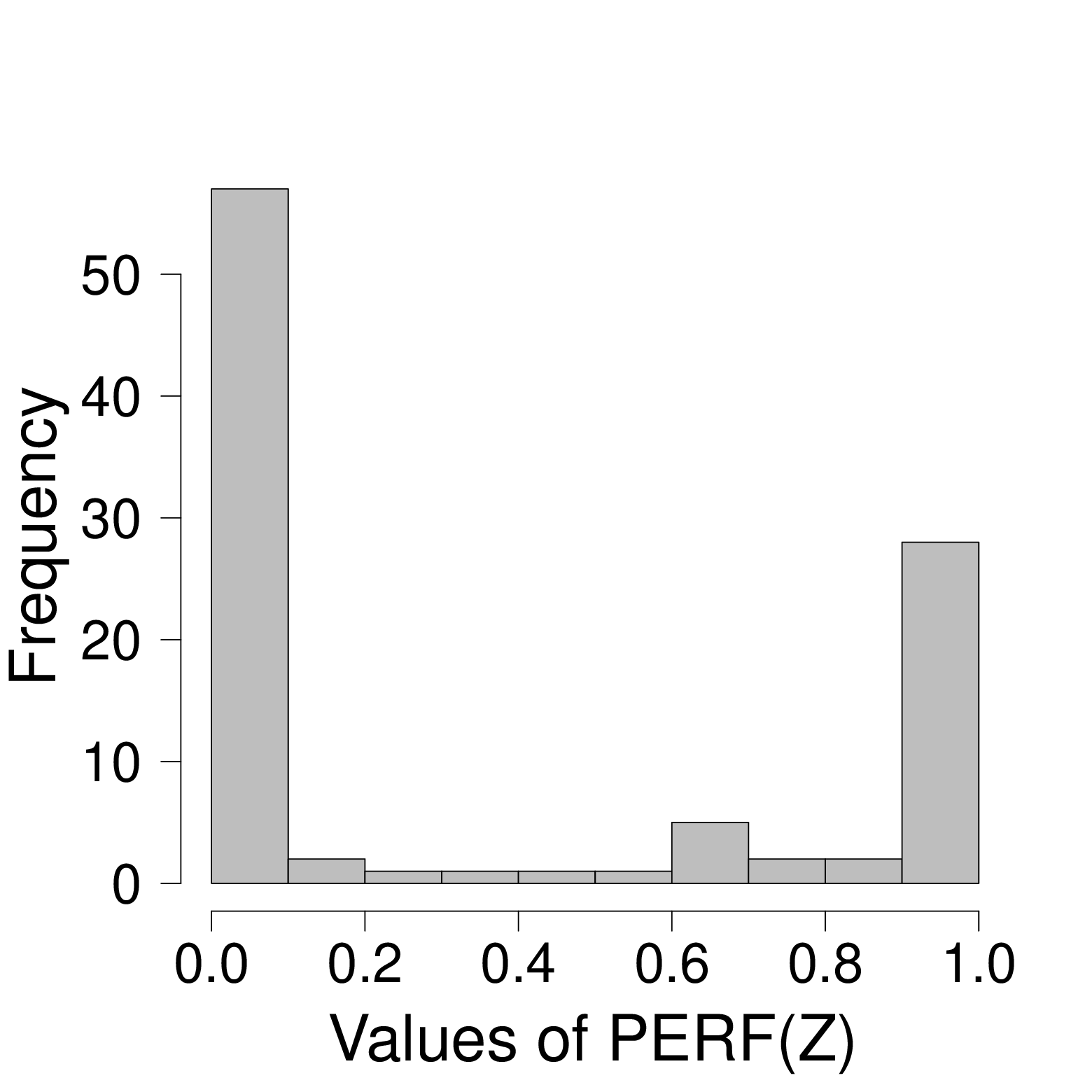}
\includegraphics[width=1.5in]{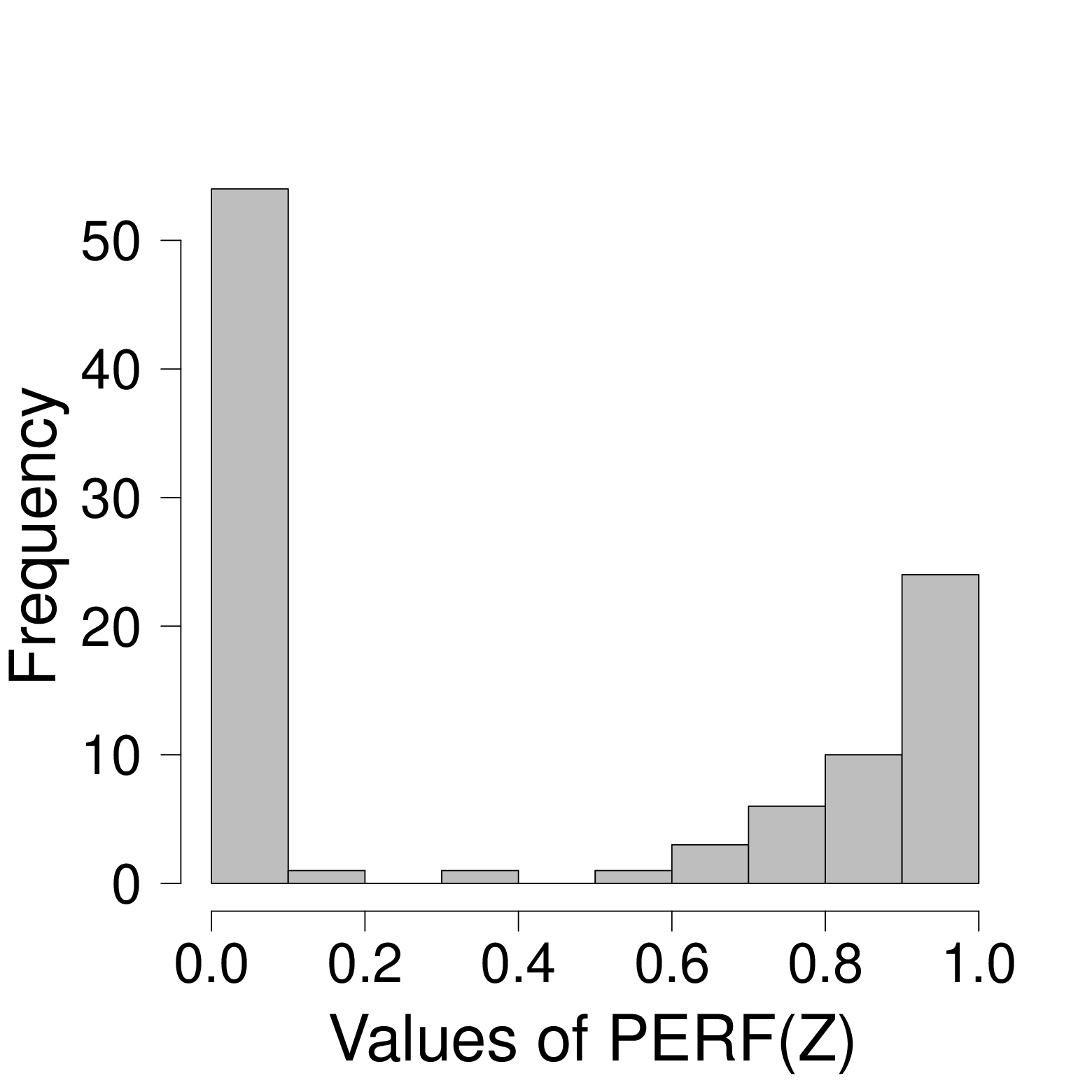}
\includegraphics[width=1.5in]{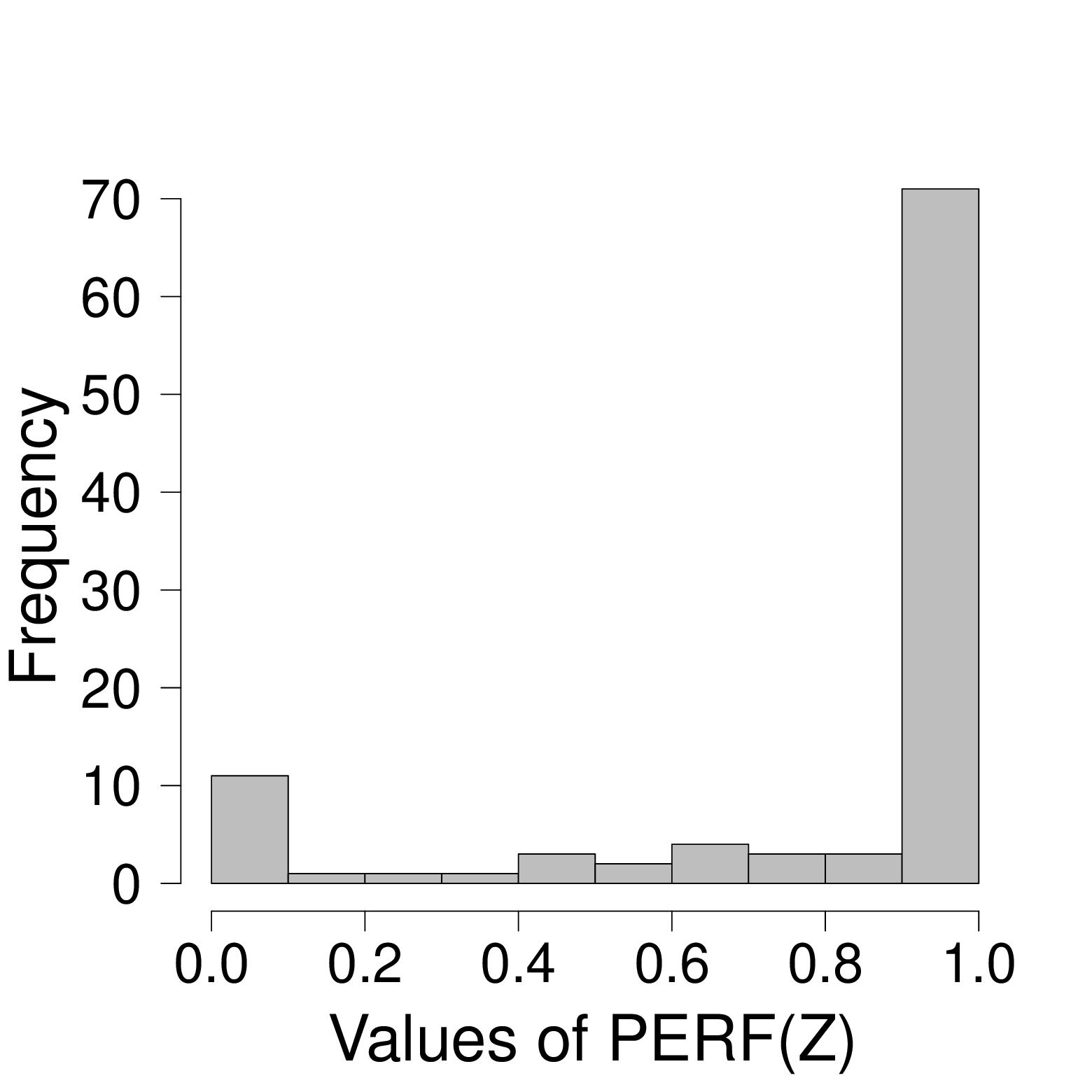}
\includegraphics[width=1.5in]{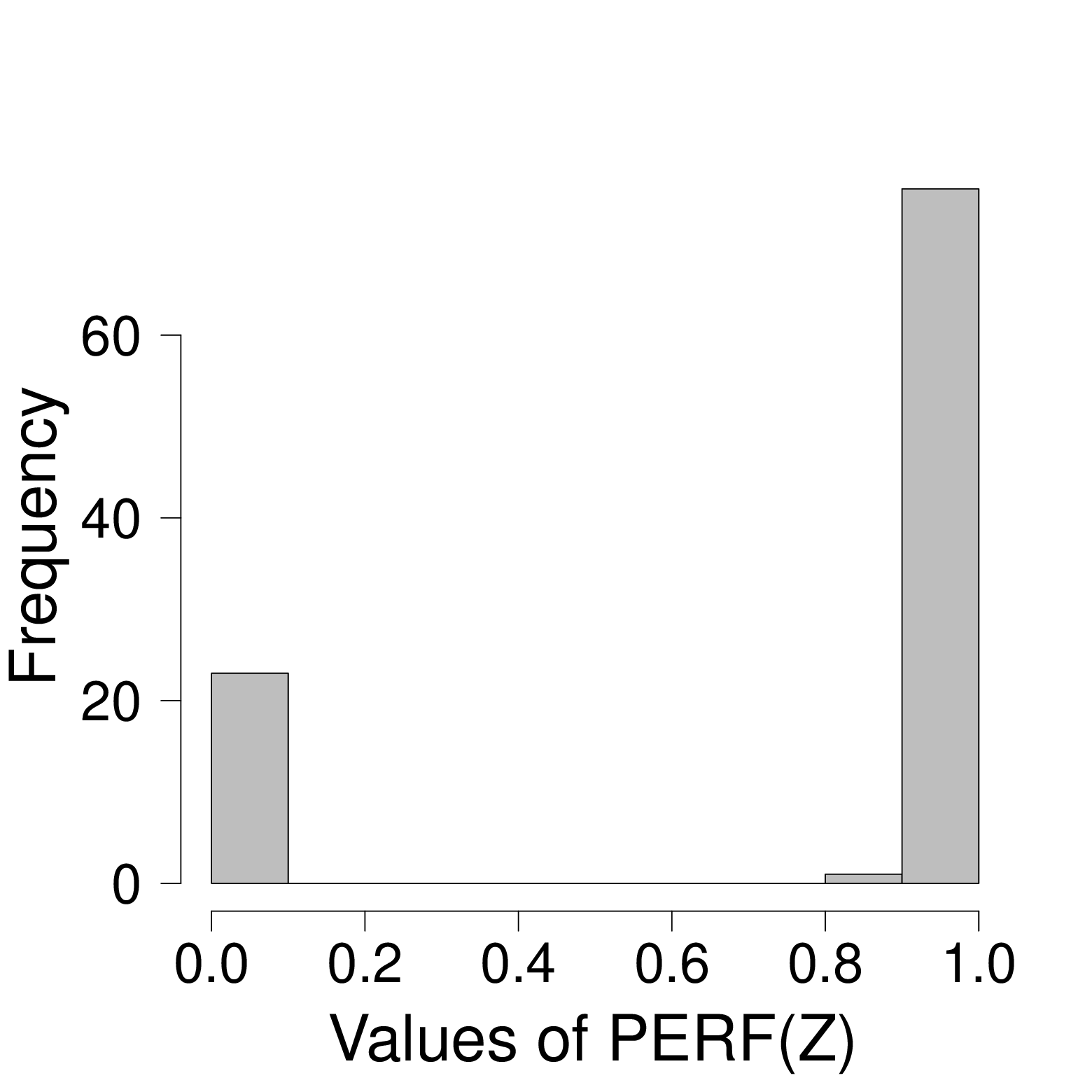}
\\
\includegraphics[width=1.5in] {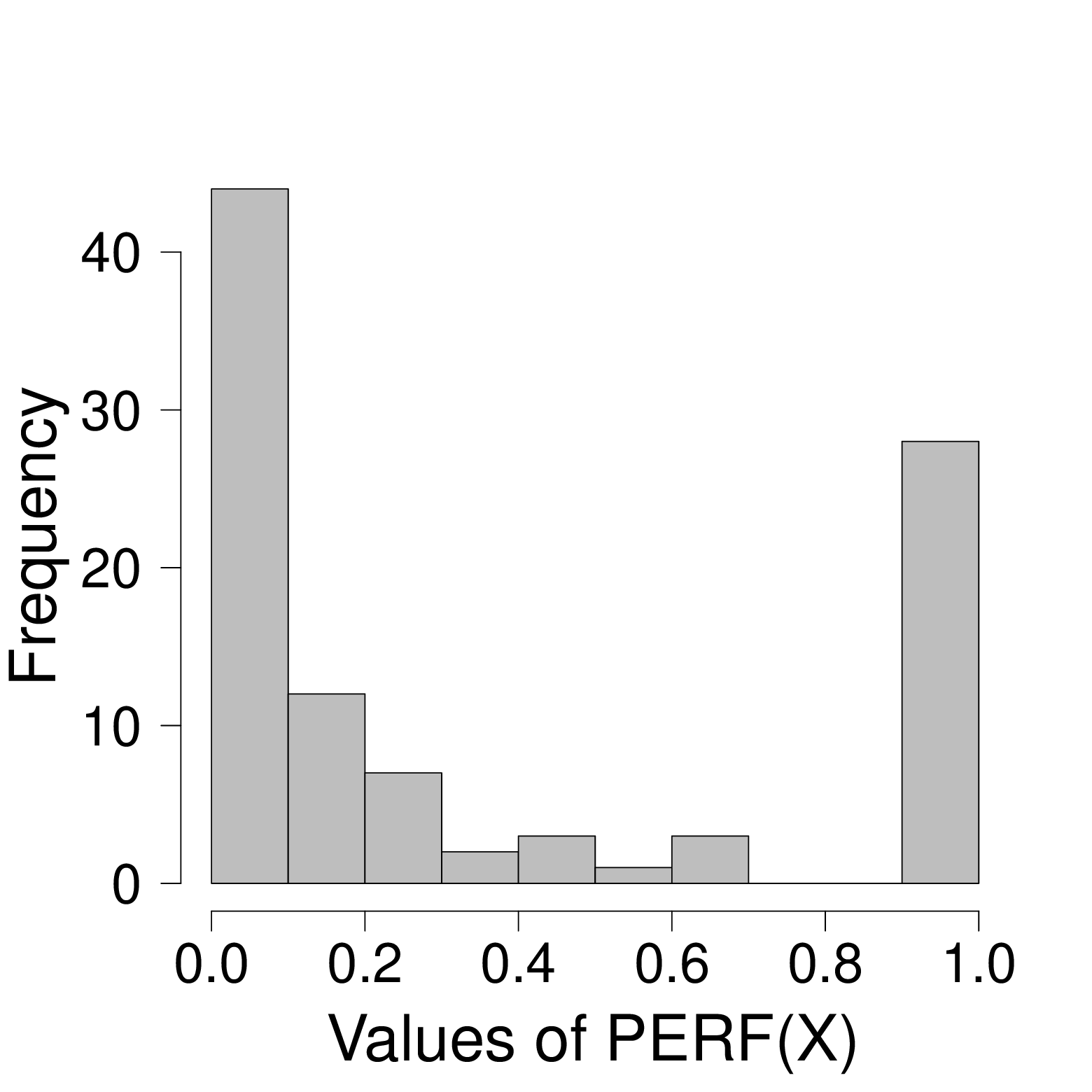}
\includegraphics[width=1.5in] {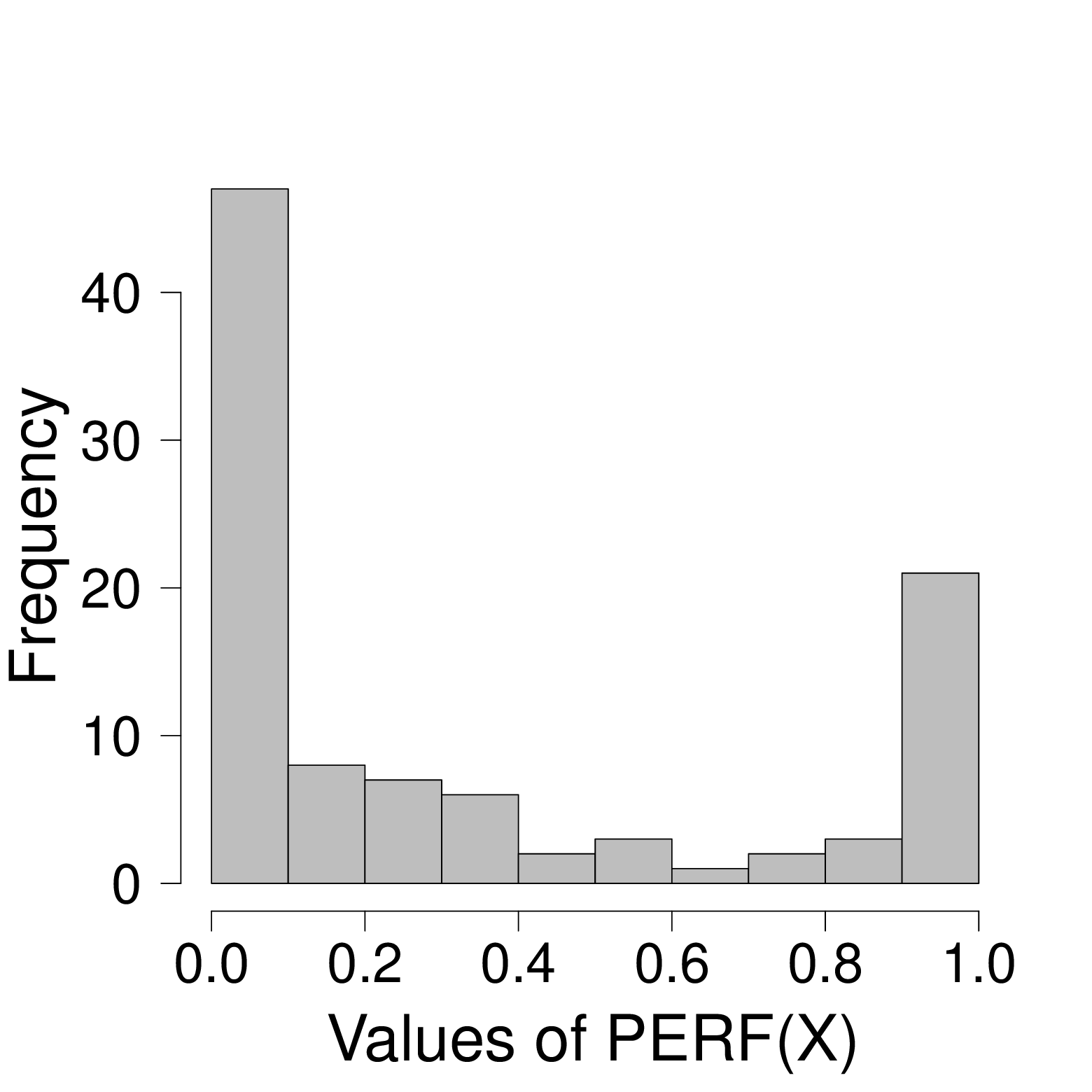}
\includegraphics[width=1.5in] {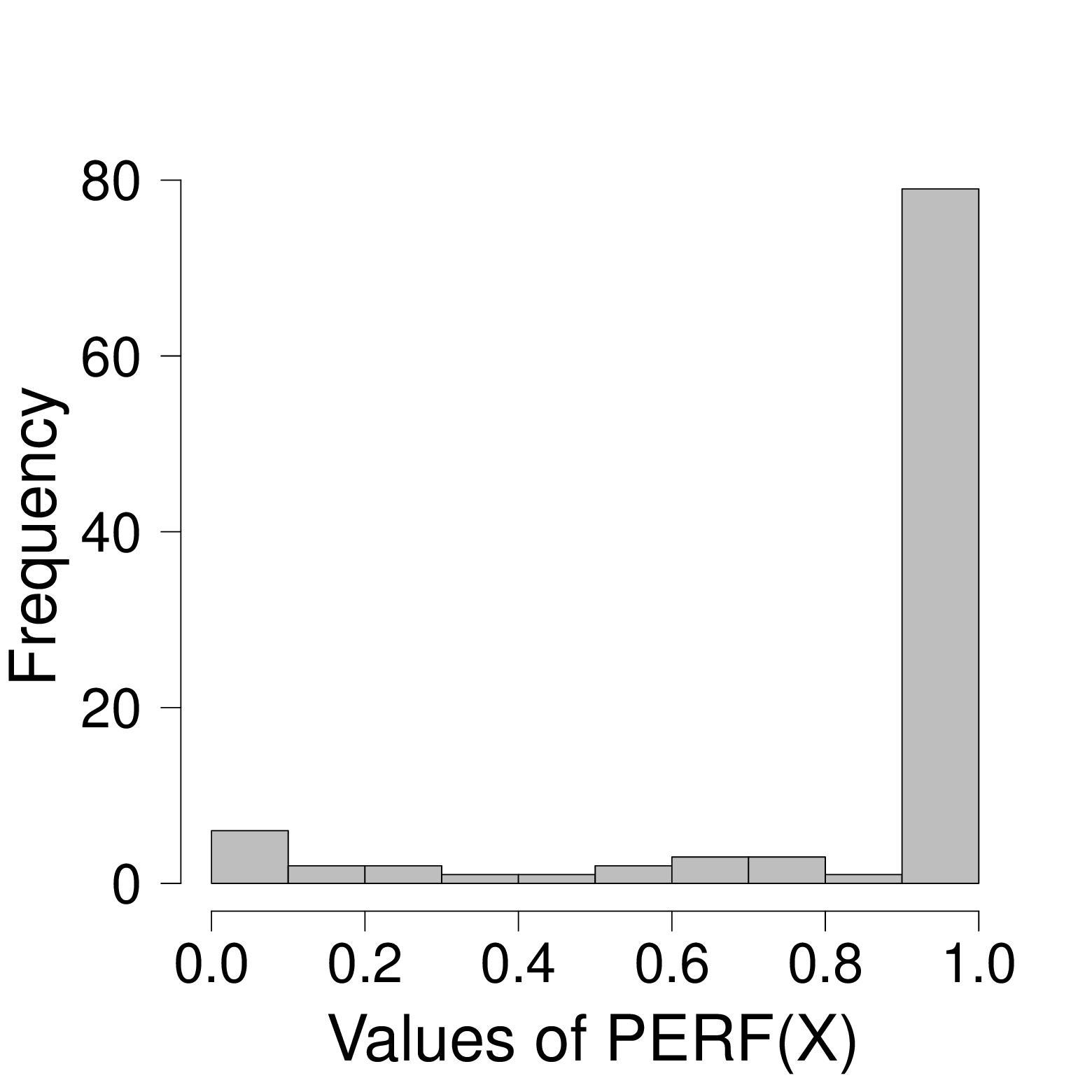}
\includegraphics[width=1.5in] {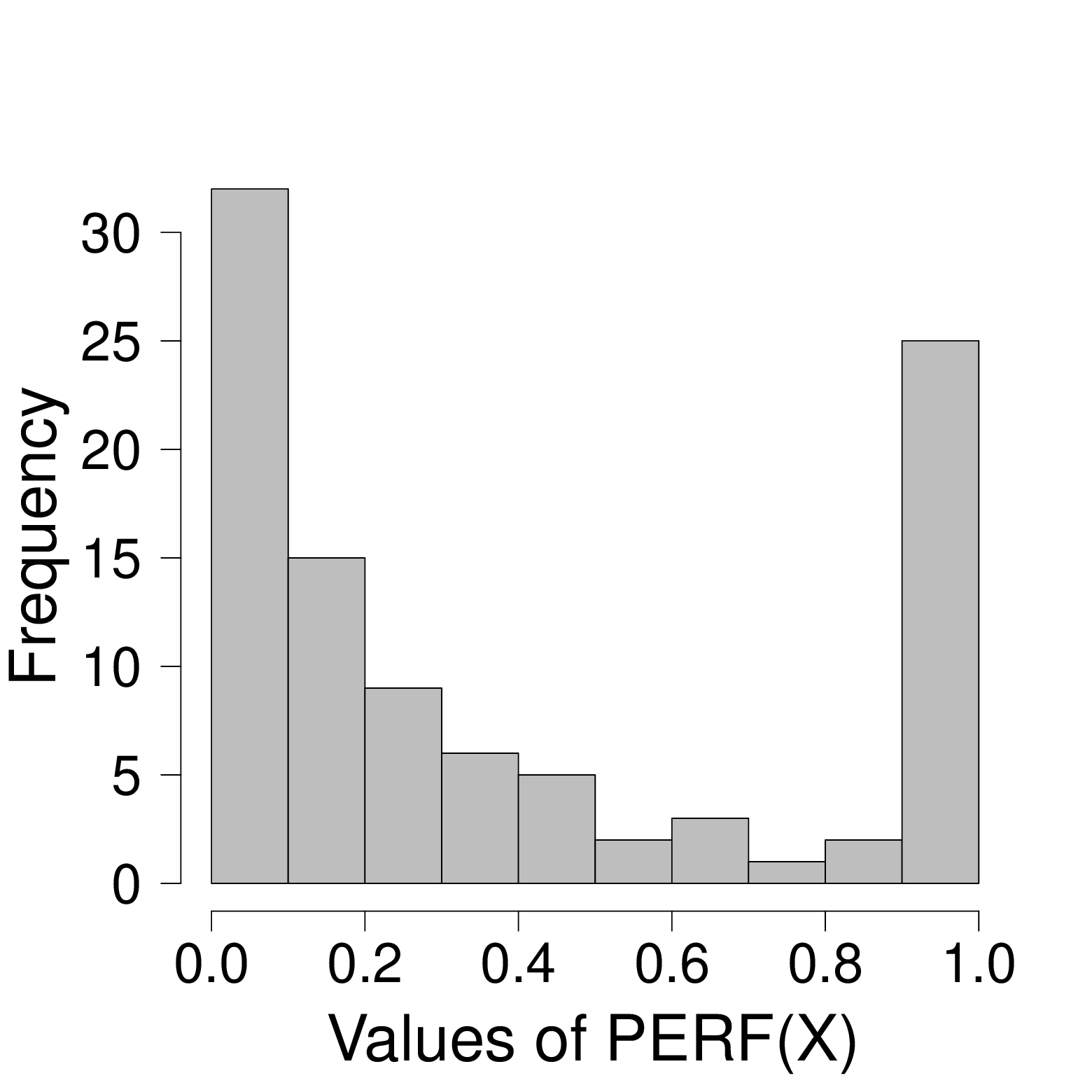}
\vglue-40pt
~
\begin{center}
\begin{minipage}[t]{0.90\textwidth}
\caption{\footnotesize Study of the distribution of $PERF$ when $p=1000$ and $s=1$. Top: study of $PERF(Z)$; Bottom: study of $PERF(X)$. From left to right: $(n,m)  = (5,20)$, $\rho = 0.1$ and $\sigma^2 = 1$; $(n,m)  = (5,20)$, $\rho = 0.9$ and $\sigma^2 = 1$; $(n,m)  = (5,20)$, $\rho = 0.9$ and $\sigma^2 = 25$; $(n,m)  = (5,500)$, $\rho = 0.9$ and $\sigma^2 = 1$.}
\label{fig:histMultiplp1000s1}
\end{minipage}
\end{center}
\end{figure}

\noindent{\it The high dimensional case.} Table~\ref{tab:ResultGdDim}, Figure~\ref{fig:histMultiplp50s20} and Figure~\ref{fig:histMultiplp1000s1} summarize the results in this case.
The main observation is that the behavior of the quantities $PERF(Z)$, $PERF(X)$ and $PERF(I)$ highly depends on whether the sparsity index $s$ is larger than the sample size $n$ or not.
Let us distinguish these two cases.\\
\noindent $-$ {\it When $s>n$:} in this difficult case, the performance of the Transductive LASSO varies with the dimension $p$.
Indeed, for moderates $p$ (a dimension smaller than about $100$), the Transductive LASSO has good performance compared to the LASSO, as observed in the first three lines of Table~\ref{tab:ResultGdDim}, where $n=5<8=s$. 
Note that in this case, the improvements using the Transductive LASSO are particularly observed for the denoising error $PERF(X)$, with a median value equal to $0.23$ and $0.42$ respectively when $\rho=0.1$ and $\rho=0.9$ (and with $\sigma^2 = 1$).
Nevertheless, we notice that the transductive error ratio $PERF(Z)$ is altered when the noise level increases.
In this case the denoising error ratio $PERF(X)$ is even more affected (see line~3 in Table~\ref{tab:ResultGdDim}).
Despite this alteration, the Transductive LASSO has still a nice behavior in this setting. The same conclusions can be made in the experiments related to Figure~\ref{fig:histMultiplp50s20}, where $p=50$, $s=20$ and $n=10$. Note that in this example, and when the parameters $\rho=0.1$ and $\sigma^2=1$, the median values of the $PERF(Z)$ and $PERF(X)$ are respectively $0.82$ and $0.14$.\\
\noindent On the other hand for large $p$ (and when $s>n$), it turns out that the Transductive LASSO does improve the LASSO estimator only when the variables are correlated. This is illustrated in Table~\ref{tab:ResultGdDim} (line~4 and~5) when $p=1000$, $s=50$ and $n=20$.

\noindent $-$ {\it When $s\leq n$:} we can consider two sub-cases, making a difference between problems with a large sparsity index $s$ (in comparison to $p$) and the others with a small one.
The last three columns of Table~\ref{tab:ResultGdDim} summarize the performance of the methods when $s$ is large.
First, note that the Transductive LASSO improves poorly the LASSO in terms of the transductive error $PERF(Z)$ when $\rho=0.1$ and $\sigma^2=1$.
Nevertheless, it remains satisfying when we deal with the prediction error $PERF(X)$.
On top of that, the Transductive LASSO, seems to be particularly interesting when the predictors are highly correlated (line~7 in Table~\ref{tab:ResultGdDim} with $\rho = 0.9$), even in presence of noise (last line in Table~\ref{tab:ResultGdDim} where $\rho=0.9$ and $\sigma^2=25$).
In this case, increasing the correlations between variables $\rho$ and the noise level $\sigma^2$ seems to imply better performance of the Transductive LASSO compared to the LASSO.
On the other hand, Figure~\ref{fig:histMultiplp1000s1} illustrates the case where the sparsity index $s$ is small compared to $p$.
Here, $p=1000$ and $s=1$.
It turns out that the Transductive LASSO is either very useful, or useless.
Indeed, as observed in the displayed histograms, the distribution of the quantities $PERF(Z)$ (top) and $PERF(X)$ (bottom) are mainly concentrated around $0$ (meaning very big improvement using the Transductive LASSO) and around $1$ (meaning almost no improvement using the Transductive LASSO).
The Transductive LASSO significantly improves the LASSO in general. Nevertheless the degradation of the behavior of the Transductive LASSO is here sensitive to the increase of $\sigma^2$.
One can compare for this purpose the third column in Figure~\ref{fig:histMultiplp1000s1} and the last line of Table~\ref{tab:ResultGdDim}.\\
In the high dimensional setting, increasing the size $m$ of the unlabeled dataset is not advantageous to the performance of the Transductive LASSO in terms of the transductive error.
This can be observed in the last column of Figure~\ref{fig:histMultiplp1000s1}.

\vspace{0.8cm}

\noindent {\bf Conclusion of the simulation study:}
the Transductive LASSO seems to be a good alternative to the LASSO in most of the cases.
It responds a good way not only to the Transductive objective (through $PERF(Z)$), but also to the denoising and the estimation ones (through $PERF(X)$ and $PERF(I)$ respectively).
The Transductive LASSO is particularly useful in the difficult situation, that is when the variables are highly correlated.
It is also often robust while varying the noise level.
Moreover, it appears that in general, a large amount of unlabeled dataset $m$ does not help to make the Transductive LASSO better than the LASSO.
The methods works better with small values of $m$.
Hence it turns out that more clever ways to exploit the unlabeled points can be imagined.
For instance, one may add weights to the observations.
More precisely, one can associate to each labeled point a weight, bigger than the weight set for the unlabeled points.
This would be the topic of a future work.
Furthermore, the simulation study reveals how beneficial can be the use of the unlabeled points even to increase the performance in the denoising task.\\
Finally a surprising observation in most of our experiments is that as often as not, the minimum in
$$
\min_{(\lambda_{1},\lambda_{2})\in\Lambda^{2}}
\| Z(\hat{\beta}^{TL}(\lambda_{1},\lambda_{2})-\beta^*) \|_{2}^{2}
<
\min_{(\lambda_{1},0)\in\Lambda^{2}}
\| Z(\hat{\beta}^{TL}(\lambda_{1},0)-\beta^*) \|_{2}^{2} ,
$$
{\it does not significantly depend on} $\lambda_{1}$ for a very large range of
values $\lambda_{1}$.
This is quite interesting for a practitioner as it means that in the use of the Transductive LASSO, we can reduce significantly the computation cost and deal (almost) with only a singular unknown tuning parameter (that is $\lambda_{2}$) rather than with two.

\vspace{0.8cm}

\begin{figure}[ht]
\begin{center}
\includegraphics[width=7cm, height=4.5cm]{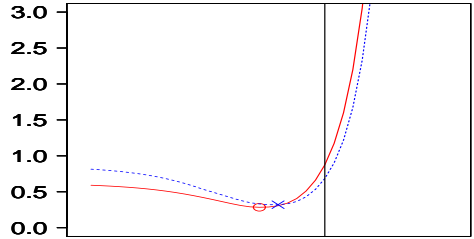}
\caption{\label{fig:perf} Evolution of the denoising error (the red solide line) and the transduction error (the blue dashed line) of the LASSO w.r.t. $\lambda$. The minimum of the denoising and the transduction errors are marked respectively by a red circle and a blue cross. The best tuning parameter for the variable selection purpose is pointed by a vertical line.}
\end{center}
\end{figure}

\noindent{\bf Discussion on the regularization parameter.}
We would like to point out the importance of the tuning parameter $\lambda$ in a general term.
Figure~\ref{fig:perf} illustrates a graph of a typical experiment in the low dimensional setting.
There are two curves on this graph, that represent the quantities $(1/n)\| X(\hat{\beta}_{\lambda}^{L}-\beta^*)\|_{2}^{2}$ and $(1/m)\| Z(\hat{\beta}_{\lambda}^{L}-\beta^*)\|_{2}^{2}$ with respect to $\lambda$.
We observe that both functions do not
reach their minimum value for the same value of $\lambda$ (the minimum are highlighted on the graph by a circle and a cross), even if these
minimum are quite close.
Since we consider variable selection methods, the identification of the true support
$\{j:\,\beta_j^*\neq 0\}$ of the vector $\beta^*$ is also in concern. One expects that the
estimator $\hat{\beta}$ and the true vector $\beta^*$ share the same support at least when
$n$ is large enough. This is known as the variable selection consistency problem and it
has been considered for the LASSO estimator in several works \cite{Bunea_consist,MeinshBulhmConsistLasso,MeinYuSelect,WainSelection,BiYuConsistLasso}.
Recently, \cite{KarimNormSup} provided the variable selection consistency of the Dantzig
Selector. Other popular selection procedures, based on the LASSO estimator, such as the
Adaptive LASSO~\cite{AdapLassoZou}, the SCAD~\cite{FanLiScad}, the S-LASSO~\cite{Mo7SLasso}
and the Group-LASSO~\cite{BachGpLasso}, have also been studied under a variable selection point of view.
Following our previous work \cite{L1MOH}, it is possible to provide such results for
the Transductive LASSO.
The variable selection task has also been illustrated in Figure~\ref{fig:perf} by the vertical line.
We reported the minimal value of $\lambda$ for which
the LASSO estimator identifies correctly the non zero components of $\beta^*$. This value of
$\lambda$ is quite different from the values that minimizes the prediction loss.
This observation is recurrent
in almost all the experiments: the estimation $X\beta^*$, $Z\beta^*$ and the support of $\beta^*$ are three
different objectives and have to be treated separately. We cannot expect in general to find a choice for $\lambda$ which makes the LASSO, for instance, has good performance for all the mentioned objective simultaneously.

\subsection{Real data}
\label{sec:RealData}

\begin{table}[t]
\caption{Evaluation of the the median and the quantile of order $0.3$ ($Med[Q_3]$) of the quantities $PERF(Z)$ in the high dimensional real dataset. Here $n$ is the labeled sample size and $t = m-n$ is the unlabeled sample size.}
\label{tab:RealDataPERFZ}
\begin{center}
\begin{sc}
\begin{tabular}{|l||c|c|c|c|c|c|}
\hline
\backslashbox{ $n$ }{ $t$ } & $10$ & $20$  & $50$   & $100$ & $500$ & $1000$ \\
\hline
\hline 
$10$     & $0.85\,[0.67]$ & $0.88\,[0.77]$  & $0.90\,[0.84]$   & $0.97\,[0.94]$ & $0.98\,[0.97]$  & $0.99\,[0.98]$ \\
\hline
$20$   & $0.74\,[0.52]$ &  $0.85\,[0.70]$  & $0.86\,[0.76]$  & $0.91\,[0.86]$ & $0.96\,[0.95]$  &  $0.98\,[0.97]$  \\
\hline
$50$  &   $0.78\,[0.49]$  & $0.68\,[0.53]$  & $0.81\,[0.67]$ & $0.84\,[0.72]$  &  $0.94\,[0.91]$ &  $ 0.97\,[0.96]$  \\
\hline
$100$ & $0.72\,[0.47]$  &  $0.68\,[0.47]$ & $0.75\,[0.58]$    & $0.75\,[0.63]$  &  $0.87\,[0,84]$  & $0.95\,[0.93]$ \\
\hline
$500$  &   $0.43\,[0.29]$  & $0.51\,[0.34]$  & $0.49\,[0.30]$ & $0.49\,[0.39]$  &  $0.81\,[0.73]$ &  $ 0.88\,[0.83]$  \\
\hline
\hline
$1000$ & $0.96\,[0.93]$  &  $0.97\,[0.92]$ & $0.91\,[0.85]$    & $0.89\,[0.81]$  &  $0.88\,[0.83]$  & $0.86\,[0.80]$ \\
\hline
\end{tabular}
\end{sc}
\end{center}
\vskip -0.1in
\end{table}

\begin{table}[t]
\caption{Evaluation of the the median and the quantile of order $0.3$ ($Med\,[Q_3]$) of the quantities $PERF(X)$ in the high dimensional real dataset. Here $n$ is the labeled sample size and $t = m-n$ is the unlabeled sample size.}
\label{tab:RealDataPERFX}
\begin{center}
\begin{small}
\begin{sc}
\begin{tabular}{|l||c|c|c|c|c|c|}
\hline
\backslashbox{ $n$ }{ $t$ } & $10$ & $20$  & $50$   & $100$ & $500$ & $1000$ \\
\hline
\hline 
$10$     & $0.04\,[0.01]$ & $0.010\,[0.003]$  & $0.04\,[0.02]$   & $0.005\,[0.001]$ & $0.0011\,[0.0003]$  & $0.12\,[0.10]$ \\
\hline
$20$   & $0.030\,[0.004]$ &  $0.006\,[0.002]$  & $0.004\,[0.001]$  & $0.005\,[0.002]$ & $0.003\,[0.001]$  &  $0.13\,[0.11]$  \\
\hline
$50$  &   $0.07\,[0.01]$  & $0.028\,[0.008]$  & $0.012\,[0.002]$ & $0.011\,[0.003]$  &  $0.014\,[0.005]$ &  $ 0.15\,[0.13]$  \\
\hline
$100$ & $0.32\,[0.41]$  &  $0.07\,[0.01]$ & $0.024\,[0.005]$    & $0.029\,[0.009]$  &  $0.02\,[0.01]$  & $0.22\,[0.18]$ \\
\hline
$500$  &   $0.40\,[0.22]$  & $0.44\,[0.25]$  & $0.33\,[0.11]$ & $0.24\,[0.08]$  &  $0.38\,[0.30]$ &  $ 0.51\,[0.46]$  \\
\hline
\hline
$1000$ & $0.97\,[0.93]$  &  $0.97\,[0.94]$ & $0.92\,[0.86]$    & $0.90\,[0.80]$  &  $0.87\,[0.77]$  & $0.78\,[0.67]$ \\
\hline
\end{tabular}
\end{sc}
\end{small}
\end{center}
\vskip -0.1in
\end{table}

\begin{table}[t]
\caption{Evaluation of the mean ($Mean$), the median ($Med$) and the quantile $Q_3$ of order $0.3$ of the quantities $PERF(Z)$, $PERF(X)$, when the methods are used in the high dimensional real dataset.}
\label{tab:RealData}
\begin{center}
\begin{footnotesize}
\begin{sc}
\begin{tabular}{|l||c|c|c||c|c|r|}
\cline{2-7}
 \multicolumn{1}{c|}{}  & \multicolumn{3}{c||}{$PERF(Z)$}    & \multicolumn{3}{c|}{$PERF(X)$}     \\
\hline Gene     & $Moy$ & $Med$  & $Q_3$   & $Moy$ & $Med$  & $Q_3$ \\
\hline
\hline
$1$ &  $0.91$  & $0.94$  & $0.88$ & $0.44$  &  $0.39$  & $ 0.25$ \\
\hline
$2 $  &   $0.90$  & $0.92$  & $0.87$ & $0.41$  &  $0.42$ &  $ 0.21$  \\
\hline
Random & $0.88$  &  $0.92$ & $0.83$    & $0.34$  &  $0.27$  & $0.09$ \\
\hline
\end{tabular}
\end{sc}
\end{footnotesize}
\end{center}
\vskip -0.1in
\end{table}

We apply the Transductive LASSO and the LASSO estimators to a genetic dataset, where the goal is to learn the complex combinatorial code underlying gene expression.
These data have already been analyzed in~\cite{MB07} and the original source is~\cite{BT04}.
The problem we consider here is known as motif regression~\cite{CLLL03}.
By motif, we think of a sequence of letters consisting of A, C, G and T.
The instances in this dataset are genes coming from yeast.
More precisely, $L = 2587$ genes are available.
Also we have $p=666$ variables.
Each of them (with length $2587$) consists of scores associated to a given candidate motif and are computed.
These scores measure how well the motifs are represented in the upstream regions of the genes.
To summary, each row of this $L\times p$ design matrix corresponds to a gene and each column to a motif score.
In other words, each component $(i,j)\in\{1,\ldots,L\}\times\{1,\ldots,p\}$ of this matrix measures how well the $j$-th motif score is represented in the upstream region of the $i$-th gene.
The response vector is a vector of size $L$.
Its $i$-th component is the expression value of the $i$-th gene.
Actually, $255$ response vectors are available.
These several measurements have been collected based on a time-course experiment.
Then, each response vector corresponds to a measurement of the gene expressions at a time-point.
In our study, we use only one response vector by experiment.
Then we first pick one of the $255$ time-points.
According to the construction of the labeled and the unlabeled datasets, we choose to pick each of them randomly among the $2587$ available instances.\\
In the first experiment, we only consider the vector corresponding to the first time-point.
Then, we construct $X$, $Y$ and $Z$.
We first pick $n$ observations with the corresponding labels to construct $X$ and $Y$ respectively.
In order to build $Z$, we add $t= m-n$ other observations (for which we do not care about the corresponding labels) to $X$.
The values of $n$ and $t$ are specified in Tables~\ref{tab:RealDataPERFZ} (for transductive the error) and~\ref{tab:RealDataPERFX} (for the denoising error), where the results for this setting are summarized.\\
\noindent Most of these results confirm what has been observed in the simulation study.
Indeed, we remark a difference in the performance of the methods in the high dimensional case and when $p<n$ (we recall that $p=666$).
The difference between the last line, where $n=1000$, and the other lines of both Tables~\ref{tab:RealDataPERFZ} and~\ref{tab:RealDataPERFX} illustrates this point.
Indeed, when $n$ is large, the improvement using the Transductive LASSO is not that significant for both the transductive and the denoising errors (about $0.90$).
We observe a big difference with the high dimensional case (the lines above), where the improvement using the Transductive LASSO is to be noticed most of the time.
Conforming to the simulation study, the performance of the Transductive LASSO are particularly marked for the denoising error.
Indeed, $PERF(X)$ is very low, with a median value between $0.001$ and $0.50$, as displayed in Table~\ref{tab:RealDataPERFX}.
Moreover, the performance of the Transductive LASSO compared to the LASSO are getting better and better when $n$ is small.
According to the transductive error (Tables~\ref{tab:RealDataPERFZ}), we also observe that the Transductive LASSO improves the LASSO estimator.
Also conforming to the simulation study, it turns out that the improvement using the Transductive LASSO is not that significant when $t$ (and then $m$) is large.
Actually, the best case in this real dataset corresponds to the situation where $n$ is large ($n=500$) and $t$ is small ($t=10$), with a median value of $PERF(Z)$ equal to $0.43$.
Another observation can be made.
According to the results displayed in Tables~\ref{tab:RealDataPERFZ}, we remark the diagonal (with $n=t$) plays an important role.
Indeed, the value of $PERF(Z)$ when $n=t$ is around $0.8$.
Moreover, when $n>t$ the improvement is always better than $0.8$ in these high dimensional experiments.
This let us believe that the best situations for the Transductive LASSO here, but also in general, is when $n>t$.\\
In all these results, we expect that the sparsity index $s$ played a role.
Indeed, we already have seen in the simulation experiments that the cases where $n>s$ and those where $n<s$ are different.
Nevertheless, our above study does not able us to make a conclusion on an approached value of $s$.

\vspace{0.2cm}

Let us now consider the second study.
Here the way to construct $X,\,Y$ and $Z$ is the same as previously, excepted for the the values of $n$ and $m$.
Here both of them are random in $1,\ldots,L$ (recall that $L=2587$ is the total number of the available instances) and such that $L\geq m >2n$.
Then it is the less advantageous situation for the Transductive LASSO.
These results can be then associated to the upper diagonal results of Tables~\ref{tab:RealDataPERFZ} and~\ref{tab:RealDataPERFX}.
The main aspect of this study is that the time-point differs.
Indeed, we choose the first time-point in the first experiment, the second in the second study, whereas we pick randomly one time-point for each replication in the third experiment (cf. Table~\ref{tab:RealData}).\\
\noindent The results are summarized in Table~\ref{tab:RealData}.
This study reveals that the behavior of the Transductive LASSO compared to the LASSO remains the same for all the time-points.
We observe that even in this real dataset, the Transductive LASSO is useful.
Moreover, as expected in this case, the Transductive LASSO outperforms the LASSO estimator particularly in terms of the prediction error.

\section{Theoretical results}
\label{thms}

In this section, we consider the theoretical properties of the Transductive LASSO and Transductive Dantzig Selector, and more generally of the estimator $\hat{\beta}_{A,\lambda}$ and $\tilde{\beta}_{A,\lambda}$ given respectively by~\eqref{dfntlasso} and~\eqref{dfntds} for any given matrix $A$ ($A=\sqrt{n/m}Z$ is then a special case).

\subsection{Assumptions}

Here, we give our two assumptions. The first one is about the matrix $A$, the
second one is about the preliminary estimator $\widehat{A\beta^{*}}$.
\begin{description}
	\item[Assumption $H(A,\tau)$:] there exists a constant $c(A,\tau)>0$ such that, for any
	$\alpha\in\mathds{R}^{p}$ such that
		$ \sum_{j:\beta^{*}_{j} = 0}
		\left|\alpha_{j}\right| \leq \tau \sum_{j:\beta^{*}_{j} \neq 0}
		 \left|\alpha_{j}\right|, $
	we have
	\begin{equation}
	\label{eq:CondHypEase}
		\alpha'(A'A) \alpha
		\geq
		c(A,\tau) n \sum_{j:\beta_{j}^*\neq 0}\alpha_{j}^{2} .
	\end{equation}
\end{description}
First, let us explain briefly the meaning of this hypothesis.
In the case where $A$ has full rank, the condition
$$
\alpha'A'A \alpha
\geq
c(A,\tau) n \sum_{j:\beta_{j}^*\neq 0} \alpha_{j}^{2},
$$
is always satisfied for any $\alpha \in\mathds{R}^{p}$ with $c(A,\tau)$ larger
than the smallest eigenvalue of $A'A/n$.
However, for the LASSO, we have $(A'A)=(X'X)$ and $A'A$ cannot be invertible if $p>n$.
Even in this high dimensional setting, Assumption $H(A,\tau)$ may still be satisfied.
Indeed, the assumption requires that Inequality~\eqref{eq:CondHypEase} holds only for a small for a small subset of $\mathds{R}^{p}$ determined by the condition $ \sum_{j:\beta^{*}_{j} = 0}
          \left|\alpha_{j}\right| \leq \tau \sum_{j:\beta^{*}_{j} \neq 0}
                    \left|\alpha_{j}\right| .$
For $A=X$, this assumption becomes exactly the one taken in~\cite{Lasso3}.
In that paper, the necessity of such an hypothesis is also discussed.
\begin{description}
	\item[Assumption~${\rm conf}(\widehat{A\beta^{*}},\kappa,\eta)$:]
	The estimator $\widehat{A\beta^{*}}$ is such that, with probability at least $1-\eta$,
	$$
	\left\|A'(\widehat{A\beta^{*}}-A\beta)\right\|_{\infty}
	\leq
	\kappa \sigma
	\sqrt{2n\log\frac{p }{\eta}}.
	$$
\end{description}
This assumption will be discussed for different types of preliminary estimators
in Section~\ref{exam}.
However note that it always holds when $A=X$ and $\widehat{A\beta^{*}}=Y$ (that is, in the "usual" LASSO case).
The idea of such an assumption results from the geometrical
considerations in our previous work on confidence regions \cite{CSEL,L1MOH}.
It just means that the preliminary estimator $\widehat{A\beta^{*}}$ may be used to build a suitable confidence region for $A\beta^{*}$.

\subsection{Main results}
\label{Sec:ResultTh}

First, Theorem~\ref{thmds} below states that the estimator $\tilde{\beta}_{A,\lambda}$ satisfies a Sparsity Inequality with high probability. A particular consequence of this result is the fact that the Transductive Dantzig Selector $\tilde{\beta}_{\sqrt{\frac{n}{m}}Z,\lambda}$ satisfies a similar SI and responds to the transductive objective.
\begin{thm}
\label{thmds}
	Let us assume that Assumption $H(A,1)$ and Assumption ${\rm conf}(\widehat{A\beta^{*}},\kappa,\eta)$ are satisfied.
	Let us choose
	$$
	\lambda
	=
	\kappa \sigma \sqrt{2n\log\frac{p }{\eta}},
	$$
	for some $\eta\in]0,1[$.
	Then, with probability at least $1-\eta$, we have simultaneously
	$$
	\left\|A\left(\tilde{\beta}_{A,\lambda}-\beta^{*}\right)\right\|_{2}^{2}
	\leq
	\frac{8\kappa^{2}
	\sigma^{2}\|\beta^{*}\|_{0}}{c(A,1)}  \log\left(\frac{p}{\eta}\right)
	$$
	and
	$$
	\left\|\tilde{\beta}_{A,\lambda}-\beta^{*}\right\|_{1}
	\leq
	\frac{2\sqrt{2} \kappa \sigma \|\beta^{*}\|_{0}}{c(A,1)}
	\sqrt{\frac{ \log\left(p/\eta \right)}{n}}.
	$$
\end{thm}
We remind that all the proofs are postponed to Section~\ref{proofs} page~\pageref{proofs}. One can use this result to tackle the particular transductive task. This is the aim of Corollary~\ref{Cor:TrDanzigSI}.
\begin{cor}
 \label{Cor:TrDanzigSI}
	Let $\lambda$ be defined as in Theorem~\ref{thmds}.
	Under Assumption $H(\sqrt{\frac{n}{m}}Z,1)$ and Assumption ${\rm conf}(\widehat{\sqrt{\frac{n}{m}}Z\beta^{*}},\kappa,\eta)$, we have with probability $1-\eta$
	$$
	\frac{1}{m}\left\|Z\left(\tilde{\beta}_{\sqrt{n/m}Z,\lambda}-\beta^{*}\right)\right\|_{2}^{2}
	\leq
	\frac{8\kappa^{2} \sigma^{2} \|\beta^{*}\|_{0}}{c(\sqrt{n/m}Z,1)}
            \frac{ \log\left(p/\eta \right)}{n}.
$$
\end{cor}
Based on Theorem~\ref{thmds}, a proper choice of the matrix $A$ can also make us respond to the other objectives (denoising and estimation) we considered in Section~\ref{Sec:DiscussA}. Indeed, in those cases we obtain:
\begin{itemize}
	\item Under Assumption~$H(X,1)$ and Assumption~${\rm conf}(\widehat{X\beta^{*}},\kappa,\eta)$ and with probability at least $1-\eta$
	$$
	\left\|X\left(\tilde{\beta}_{X,\lambda}-\beta^{*}\right)\right\|_{2}^{2}
	\leq
	\frac{8\kappa^{2} \sigma^{2} \|\beta^{*}\|_{0}}{c(X,1)} 
		\log\left(\frac{p}{\eta} \right);
	$$
	\item Under Assumption~${\rm conf}(\widehat{\sqrt{n}I\beta^{*}},\kappa,\eta)$ and with probability at least $1-\eta$
	$$
	\left\|\tilde{\beta}_{\sqrt{n}I\,\lambda}-\beta^{*}\right\|_{2}^{2}
	\leq
	\frac{8\kappa^{2} \sigma^{2} \|\beta^{*}\|_{0}}{\sqrt{n}}
		\log\left(\frac{p}{\eta} \right).
	$$
\end{itemize}
Corollary~\ref{Cor:TrDanzigSI} and the above statements claim that each estimator perform well for the task it is designed to fulfill.
In a similar way, we finally can establish analog results for the Transductive LASSO and more generally for the estimator $\hat{\beta}_{A,\lambda}$ given by Definition~\ref{dfntlasso}.
\begin{thm}
\label{thmlasso}
	Let us assume that assumption $H(A,3)$ and Assumption~${\rm conf}(\widehat{A\beta^{*}},\kappa,\eta)$ are satisfied.
	Let us choose
	$$
	\lambda
	=
	2 \kappa \sigma \sqrt{2n\log\left(\frac{p }{\eta}\right)},
	$$
	for some $\eta\in]0,1[$.
	Then, with probability at least $1-\eta$, we have simultaneously
	$$ \|A(\hat{\beta}_{A,\lambda}-\beta^{*})\|_{2}^{2}
	\leq
	\frac{72 \sigma^{2}\kappa^{2} \|\beta^{*}\|_{0}}{c(A,3)}
		\log\left(\frac{p}{\eta}\right)
	$$
	and
	$$
	\left\|\beta^{*}  -\hat{\beta}_{A,\lambda}\right\|_{1}
	\leq
	\frac{24\sqrt{2} \|\beta^{*}\|_{0}}{c(A,3)} \sqrt{\frac{\log\left(p/\eta \right)}{n}}.
	$$
\end{thm}

\subsection{Examples of preliminary estimators}
\label{exam}

In this section, we examine some preliminary estimators $\widehat{A\beta^{*}}$ and check if they may satisfy Assumption ${\rm conf}(\widehat{A\beta^{*}},\kappa)$. This is an important issue of the paper since it helps to understand how restrictive are the assumptions in the results of Section \ref{Sec:ResultTh}. The first example deals with the (generalized) least square estimator.
\begin{thm}
\label{thmpreliminary}
	Let us choose $\widetilde{(X'X)}^{-1}$ any pseudo-inverse of $(X'X)$ and let us set
	$$
	\widehat{A\beta^{*}}
	=
	A \widetilde{(X'X)}^{-1} X' Y,
	$$
	as preliminary estimator.
	Then, under the assumption $\Ker(A)=\Ker(X)$ and for any $\eta\in]0,1[$, Assumption~${\rm conf}(\widehat{A\beta^{*}},\kappa,\eta)$ holds with $\kappa = \sqrt{ \frac{p}{\sum_{j=1}^p\widetilde\Omega_{j,j}^{-1}  }}$ where $\widetilde\Omega = ((A'A)\widetilde{(X'X)}^{-1}(A'A))/n$.
\end{thm}
According to Theorem~\ref{thmpreliminary}, the standard case of interest is when $A=X$.
The preliminary estimator becomes $\widehat{X\beta^{*}}=Y$ and we obtain that ${\rm conf}(Y,1,\eta)$ holds.
Plugging this into Theorems~\ref{thmds} and~\ref{thmlasso} implies the theorems about the LASSO and the Dantzig Selector provided in~\cite{Lasso3}.
Moreover, other choices for $A$ and $\widehat{X\beta^{*}}$ are possible which able us to deal for instance with the transductive setting.
Hence, one can interpret Theorem~\ref{thmpreliminary} together with Theorems~\ref{thmds} and~\ref{thmlasso} as a generalization of the result in~\cite{Lasso3}.\\

\noindent To introduce the second preliminary estimator, let us consider the case when $A \ne X$. Then the assumption $\Ker(A)=\Ker(X)$ is restrictive when $p>n$ (in the somehow appreciable case $p<n$, the assumption holds since both $X$ and $Z$ may have full rank).
If the relation $\Ker(A)=\Ker(X)$ is not satisfied, as the construction of $Z$ leads to $\Ker(Z)\subset\Ker(X)$, we may suggest the following alternative.
Consider the restriction of the estimation procedure to the span of $X$. That is, let replace $Z$ by $Z_{X}=\widetilde{(X'X)}^{-1}(X'X)Z$.
Then the assumption $\Ker(Z_{X})=\Ker(X)$ is satisfied.
As a consequence, with probability at least $1-\eta$, the following inequality
$$
\left\|Z_{X}\left(\tilde{\beta}_{\sqrt{n/m}Z_{X},\lambda}- \beta^{*}\right)\right\|_{2}^{2}
\leq
\frac{8\kappa^{2} \sigma^{2}\|\beta^{*}\|_{0}}{c(\sqrt{n/m}Z_{X},1)}  \log\left(\frac{p}{\eta}\right),
$$
is obtained for instance for the Transductive Dantzig Selector (an analog inequality can be written for the Transductive LASSO), under Assumption $H(\sqrt{n/m}Z_{X},1)$ and with the same choice of the tuning parameter $\lambda$ as in Theorem~\ref{thmds}.
Finally, let us remark that $(Z-Z_{X})\tilde{\beta}_{\sqrt{n/m}Z_{X},\lambda}=0$ and conclude the following result.
\begin{cor}
	Under Assumption~$H(\sqrt{n/m}Z_{X},1)$ and with the same choice of $\lambda$ as in Theorem \ref{thmds}, we have with probability at least $1-\eta$,
	$$
	\left\|Z\left(\tilde{\beta}_{\sqrt{n/m}Z_{X},\lambda}- \beta^{*}\right)\right\|_{2}^{2}
	\leq
	\frac{8\kappa^{2} \sigma^{2}\|\beta^{*}\|_{0}}{c(\sqrt{n/m}Z_{X},1)}  \log\left(\frac{p}{\eta}\right)
	+
	\|(Z-Z_{X})\beta^{*}\|_{2}^{2}.
	$$
\end{cor}
The conclusion figured out this result is quite intuitive:
when $\|(Z-Z_{X})\beta^{*}\|_{2}^{2}$ is large, the information in $X$ is not sufficient to estimate $Z\beta^{*}$.
But, if $\|(Z-Z_{X})\beta^{*}\|_{2}^{2}$ is small, the Transductive Dantzig Selector based on $Z_{X}$ has good performances.
This assumption has the same status as a regularity assumption in a non-parametric setting. Obviously, we cannot know whether $\|(Z-Z_{X})\beta^{*}\|_{2}^{2}$ is small or not.
However when it is not, it seems impossible to guaranty a good estimation.\\

\noindent The final preliminary estimator we examine here has also been studied in the experiments part (cf. Section~\ref{simu}).
Let us consider the Dantzig Selector as preliminary estimator.
Here, a quite natural assumption can be made.
It somehow says that $X'X$ and $A'A$ are not too far from each other.
\begin{thm}
\label{thmpreliminarylasso}
	Let us assume that, there is a constant $k>0$ such that for any $u\in\mathds{R}^{p}$ with $\|u\|_{1}\leq 2 \|\beta^*\|_{1}$,
	$$
	\left\|\left[(X'X)-(A'A)\right]u\right\|_{\infty}
	\leq
	k\sigma  \sqrt{2n\log(p)}.
	$$
	Let moreover $\eta\in]0,1[$ and set the preliminary estimator
	$$ \widehat{A\beta^{*}}
	=
	A \hat{\beta}^{DS}_{2\sigma\sqrt{2n\log\left(\frac{p}{\eta}\right)}} .
	$$
	Then Assumption ${\rm conf}\left(\widehat{A\beta^{*}},\kappa,\eta\right) $ is true with $\kappa=4+k$.
\end{thm}
The same result would hold as well for the LASSO as a preliminary estimator. Moreover, in this last result, one can also consider the transductive objective and consider the matrix $A = \sqrt{n/m}Z_X$ as introduced above. Such a choice helps us to provide good theoretical guaranties with very mild assumptions on the Gram matrix $X'X$.


%
%




\section{Conclusion}

In this paper, we studied transductive versions of the LASSO and the Dantzig Selector.
These new methods appeared to enjoy both theoretical and practical advantages.
Indeed, in one hand, we showed that the Transductive LASSO and Dantzig Selector satisfy sparsity inequalities with weaker assumption on the Gram matrix than the original method.
On the other hand we displayed some experimental results illustrating the superiority of the Transductive LASSO on the LASSO.
On top of that, these transductive methods are easy to compute.\\
The experimental study reveals that the Transductive LASSO is often much better than the original LASSO.
Nevertheless, when the number of unlabeled observations is much larger than the sample size, it turns out the the gain using the Transductive LASSO is reduced.
We will focus on this point in a future work.

\section{Proofs}
\label{proofs}

In this section, we give the proofs of our main results.

\subsection{Proofs of Theorems \ref{thmds} and \ref{thmlasso}}

\begin{proof}[Proof of Theorem \ref{thmds}]
First, we have obviously
\begin{multline}
\left\| A \left(\tilde{\beta}_{A,\lambda}-\beta^{*}\right)\right\|_{2}^{2}
=
\left(\tilde{\beta}_{A,\lambda}-\beta^{*}\right)'A'A
\left(\tilde{\beta}_{A,\lambda}-\beta^{*}\right)
\leq
\left\|\tilde{\beta}_{A,\lambda}-\beta^{*}\right\|_{1}
\left\|A'A
\left(\tilde{\beta}_{A,\lambda}-\beta^{*}\right)\right\|_{\infty}
\\
\leq
\left\|\tilde{\beta}_{A,\lambda}-\beta^{*}\right\|_{1}
\left\{
\left\|A'
\left(A\tilde{\beta}_{A,\lambda}-\widehat{A\beta^{*}}\right)\right\|_{\infty}
+
\left\|A'
\left(\widehat{A\beta^{*}}-A\beta^{*}\right)\right\|_{\infty}
\right\}.
\label{intermediaire1}
\end{multline}
Then, just remark that by Assumption~${\rm conf}(\widehat{A\beta^{*}},\kappa)$,
we have, with probability at least $1-\eta$,
\begin{equation}
 \label{eqpr:truebetaconstraint}
	\left\|A'(\widehat{A\beta^{*}}-A\beta^{*})\right\|_{\infty} \leq
             \kappa \sigma
                   \sqrt{2n\log\frac{p }{\eta}} .
\end{equation}
Moreover, by the definition of $\tilde{\beta}_{A,\lambda}$ (Definition \ref{dfntds}
page \pageref{dfntds}) we have
$$ \left\|A'
\left(A\tilde{\beta}_{A,\lambda}-\widehat{A\beta^{*}}\right)\right\|_{\infty}
 \leq \lambda .$$
Then, combining the fact that $\tilde{\beta}_{A,\lambda}$ minimizes $\|\cdot\|_{1}$ among all the vectors $\beta$ satisfying
$$ \left\|A'
\left(A\beta-\widehat{A\beta^{*}}\right)\right\|_{\infty}
 \leq \lambda ,$$
and the fact that thanks to~\eqref{eqpr:truebetaconstraint} and as soon as $\lambda = \kappa \sigma \sqrt{2n\log\frac{p }{\eta}}$, the vector $\beta^{*}$ satisfies the same inequality, we have
\begin{multline*}
0 \leq \|\beta^{*} \|_{1} - \|\tilde{\beta}_{A,\lambda}\|_{1}
\leq
\sum_{\beta^{*}_{j} \neq 0} |\beta^{*}_{j}| -
       \sum_{\beta^{*}_{j} \neq 0} |(\tilde{\beta}_{A,\lambda})_{j}|
         - \sum_{\beta^{*}_{j} = 0} |(\tilde{\beta}_{A,\lambda})_{j}|
\\
\leq
\sum_{\beta^{*}_{j} \neq 0} |\beta^{*}_{j}-(\tilde{\beta}_{A,\lambda})_{j}|
 - \sum_{\beta^{*}_{j} = 0} |\beta^{*}_{j}-(\tilde{\beta}_{A,\lambda})_{j}|.
\end{multline*}
As a consequence, we have
$$ \sum_{\beta^{*}_{j} = 0} |\beta^{*}_{j}-(\tilde{\beta}_{A,\lambda})_{j}|
 \leq
 \sum_{\beta^{*}_{j} \neq 0} |\beta^{*}_{j}-(\tilde{\beta}_{A,\lambda})_{j}| ,$$
which implies that the vector $\beta^{*}-\tilde{\beta}_{A,\lambda}$ is an admissible $\alpha$
for the relation in Assumption~$H(A,1)$.
Hence, using this assumption in the last above inequality, we have the following upper bound
\begin{multline}
\left\|\tilde{\beta}_{A,\lambda}-\beta^{*}\right\|_{1}
= \sum_{\beta^{*}_{j} = 0} |\beta^{*}_{j}-(\tilde{\beta}_{A,\lambda})_{j}|
 +
 \sum_{\beta^{*}_{j} \neq 0} |\beta^{*}_{j}-(\tilde{\beta}_{A,\lambda})_{j}|
\leq 2\sum_{\beta^{*}_{j} \neq 0} |\beta^{*}_{j}-(\tilde{\beta}_{A,\lambda})_{j}|
\\
\leq \sqrt{ {\rm card}\{j:\beta_{j}^{*} \neq 0\}
                \sum_{\beta^{*}_{j} \neq 0} (\beta^{*}_{j}-(\tilde{\beta}_{A,\lambda})_{j})^{2}}
\leq \sqrt{ \frac{\|\beta^{*}\|_{0}}{n c(A,1)}
    \|A
(\tilde{\beta}_{A,\lambda}-\beta^{*})\|_{2}^{2}}.
\label{intermediaire2}
\end{multline}
We plug this result into Inequality \eqref{intermediaire1} to obtain, with
probability at least $1-\eta$,
$$
\|A(\tilde{\beta}_{A,\lambda}-\beta^{*})\|_{2}^{2}
\leq
2 \kappa \sigma \sqrt{2\log\left(\frac{p }{\eta} \right)
 \frac{ \|\beta^{*}\|_{0}}{ c(A,1)  }
    \|A(\tilde{\beta}_{A,\lambda}-\beta^{*})\|_{2}^{2}},
$$
that leads to
$$
\|A(\tilde{\beta}_{A,\lambda}-\beta^{*})\|_{2}^{2}
\leq
\frac{8\kappa^{2}
\sigma^{2}}{c(A,1)} \|\beta^{*}\|_{0} \log\left(\frac{p}{\eta}\right).
$$
Plugging this last inequality into Inequality~\eqref{intermediaire2} gives
$$
\left\|\tilde{\beta}_{A,\lambda}-\beta^{*}\right\|_{1}
\leq
\frac{2\sqrt{2} \|\beta^{*}\|_{0}}{c(A,1)}
\sqrt{\frac{\log\left(\frac{p}{\eta}\right)}{n}},
$$
and this ends the proof.
\end{proof}

\begin{proof}[Proof of Theorem \ref{thmlasso}]
By the definition of the transductive LASSO (Definition \ref{dfntlasso} page
\pageref{dfntlasso}) we have
\begin{equation*}
	-2\widehat{A\beta^{*}}'A\hat{\beta}_{A,\lambda}
	+
	\hat{\beta}_{A,\lambda}'A'A\hat{\beta}_{A,\lambda}
	+
	2\lambda \|\hat{\beta}_{A,\lambda}\|_{1}
	\leq
	-2\widehat{A\beta^{*}}'A\beta^{*}
	+
	(\beta^{*})'A'A\beta^{*}
	+
	2\lambda \|\beta^{*}\|_{1}.
\end{equation*}
We can rewrite that as
\begin{multline*}
-2(\beta^{*})'A'A\hat{\beta}_{A,\lambda}
+ 2\left[A\beta^{*}-\widehat{A\beta^{*}}\right]'A\hat{\beta}_{A,\lambda}
+ \hat{\beta}_{A,\lambda}'A'A\hat{\beta}_{A,\lambda}
+2\lambda \|\hat{\beta}_{A,\lambda}\|_{1}
\\
\leq
- (\beta^{*})'A'A\beta^{*}
+ 2\left[A\beta^{*}-\widehat{A\beta^{*}}\right]'A\beta^{*}
+2\lambda \|\beta^{*}\|_{1},
\end{multline*}
or, rearranging the terms,
\begin{multline}
\label{steplasso1}
\|A(\hat{\beta}_{A,\lambda}-\beta^{*})\|_{2}^{2}
= (\hat{\beta}_{A,\lambda}-\beta^{*})A'A(\hat{\beta}_{A,\lambda}-\beta^{*})
\\
\leq
2\left[A\beta^{*}-\widehat{A\beta^{*}}\right]'A\left(\beta^{*}
-\hat{\beta}_{A,\lambda}
\right)
+2\lambda\left[\|\beta^{*}\|_{1}
  -\|\hat{\beta}_{A,\lambda}\|_{1}\right].
\end{multline}
Now, let us remark that
\begin{multline*}
\left[A\beta^{*}-\widehat{A\beta^{*}}\right]'A\left(\beta^{*}
-\hat{\beta}_{A,\lambda}
\right)
= \left[A'\left(A\beta^{*}-\widehat{A\beta^{*}}\right)\right]'\left(\beta^{*}
-\hat{\beta}_{A,\lambda}
\right)
\\
\leq
\left\|A'\left(A\beta^{*}-\widehat{A\beta^{*}}\right)\right\|_{\infty}
\left\|\beta^{*}
-\hat{\beta}_{A,\lambda}\right\|_{1}
\leq
\frac{\lambda}{2} \left\|\beta^{*}
-\hat{\beta}_{A,\lambda}\right\|_{1},
\end{multline*}
with probability $1-\eta$, provided that $\lambda=2\kappa \sigma
                   \sqrt{2n\log\frac{p }{\eta}}$ together with Assumption ${\rm conf}(\widehat{A\beta^{*}},\kappa)$.
We plug that into Inequality \eqref{steplasso1} to obtain, with probability $1-\eta$,
\begin{equation*}
\|A(\hat{\beta}_{A,\lambda}-\beta^{*})\|_{2}^{2}
\leq
\lambda \left[\left\|\beta^{*}
 -\hat{\beta}_{A,\lambda}\right\|_{1}
+ 2 \left(\|\beta^{*}\|_{1}
  -\|\hat{\beta}_{A,\lambda}\|_{1}\right)\right].
\end{equation*}
This leads to
\begin{multline}
\|A(\hat{\beta}_{A,\lambda}-\beta^{*})\|_{2}^{2}
+ \lambda \left\|\beta^{*}
 -\hat{\beta}_{A,\lambda}\right\|_{1}
\leq 2 \lambda \left(\left\|\beta^{*}
 -\hat{\beta}_{A,\lambda}\right\|_{1}
+ \|\beta^{*}\|_{1}
  -\|\hat{\beta}_{A,\lambda}\|_{1}\right)
\\
= 2\lambda \sum_{j=1}^{p} \left(
               \left|\beta^{*}_{j}
 -(\hat{\beta}_{A,\lambda})_{j}\right|+\left|\beta^{*}_{j}\right|
      - \left|(\hat{\beta}_{A,\lambda})_{j}\right|\right)
= 2\lambda \sum_{\beta^{*}_{j}\neq 0} \left(
               \left|\beta^{*}_{j}
 -(\hat{\beta}_{A,\lambda})_{j}\right|+\left|\beta^{*}_{j}\right|
      - \left|(\hat{\beta}_{A,\lambda})_{j}\right|\right)
\\
\leq 4 \lambda \sum_{\beta^{*}_{j}\neq 0} \left(
               \left|\beta^{*}_{j}
 -(\hat{\beta}_{A,\lambda})_{j}\right|\right),
\label{steplasso2}
\end{multline}
and, from this Inequality~\eqref{steplasso2}, we deduce that
$\hat{\beta}_{A,\lambda}-\beta^{*}$ is an admissible $\alpha$ vector
in Assumption~$H(A,3)$.
Then we obtain, still from \eqref{steplasso2}
and Assumption~$H(A,3)$,
\begin{multline*}
\|A(\hat{\beta}_{A,\lambda}-\beta^{*})\|_{2}^{2}
\leq
3 \lambda \sum_{\beta^{*}_{j}\neq 0} \left(
               \left|\beta^{*}_{j}
 -(\hat{\beta}_{A,\lambda})_{j}\right|\right)
\leq 6 \kappa \sigma
       \sqrt{2n\log\left(\frac{p }{\eta}\right) \|\beta^{*}\|_{0}
             \sum_{\beta_{j}^{*}\neq 0} \left(\beta^{*}_{j}
 -(\hat{\beta}_{A,\lambda})_{j}\right)^{2}}
\\
\leq 6 \kappa \sigma
       \sqrt{\frac{2\log\left(\frac{p}{\eta}\right)\|\beta^{*}\|_{0}}{c(A,3)}
            \|A(\hat{\beta}_{A,\lambda}-\beta^{*})\|_{2}^{2} }.
\end{multline*}
This last display implies
$$ \|A(\hat{\beta}_{A,\lambda}-\beta^{*})\|_{2}^{2}
   \leq \frac{72 \sigma^{2}\kappa^{2} \|\beta^{*}\|_{0}\log\left(\frac{p}{\eta}\right)}
               {c(A,3)} .$$
We plug this last result into Inequality \eqref{steplasso2} to obtain
$$
\left\|\beta^{*}
 -\hat{\beta}_{A,\lambda}\right\|_{1}
\leq
\frac{24\sqrt{2} \|\beta^{*}\|_{0}}{c(A,3)}
\sqrt{\frac{\log\left(\frac{p}{\eta}\right)}{n}}
.$$
\end{proof}

\subsection{Proofs of Theorems \ref{thmpreliminary} and \ref{thmpreliminarylasso}}

\begin{proof}[Proof of Theorem \ref{thmpreliminary}]
The proof is quite simple. As $ Y\sim\mathcal{N}(X\beta^{*},\sigma^{2} I_{n})$,
we have
$$ \widetilde{(X'X)}^{-1}X'Y - \beta^{*} \sim
       \mathcal{N}\left(0,\sigma^{2}\widetilde{(X'X)}^{-1} \right) ,$$
and so
$$ A'A\left(\widetilde{(X'X)}^{-1}X'Y - \beta^{*} \right) \sim
       \mathcal{N}\left(0, \sigma^{2}\Omega \right) ,$$
where $\Omega$ denotes the matrix $ \Omega = \Omega(A,X) = A'A \widetilde{(X'X)}^{-1} A'A .$
Let us also define
$$ V = A'A\left(\widetilde{(X'X)}^{-1}X'Y - \beta^{*} \right) .$$
Then, for any $j\in\{1,\ldots,p\}$,
$$ V_{j} \sim \mathcal{N}\left(0,\sigma^2\Omega_{j,j}\right).$$
Using a standard inequality on the tail of Gaussian variables yields
\begin{eqnarray*}
\mathbb{P}\left(| V_j | > \kappa \sigma \sqrt{2n\log (p/\eta)}\right)
\leq
\exp \left(-\frac{ \left(\kappa \sigma \sqrt{2n\log (p/\eta)} \right)^2 }{2 \sigma^2\Omega_{j,j}}\right)
=
\exp \left( - \frac{ \kappa^2  n\log (p/\eta) }{\Omega_{j,j}}\right).
\end{eqnarray*}
Then, using a union bound and the concavity of the function $x \mapsto \exp(-x)$, we easily obtain
\begin{eqnarray*}
 \mathbb{P}\left(\| V \|_{\infty} > \kappa \sigma \sqrt{2n\log (p/\eta)}\right)
& \leq &
\sum_{j=1}^p \exp \left( - \frac{ \kappa^2  n\log (p/\eta) }{\Omega_{j,j}}\right)
\\
& \leq &
p \exp \left( - \frac{1}{p}\sum_{j=1}^p \frac{ \kappa^2  n\log (p/\eta) }{\Omega_{j,j}}\right).
\end{eqnarray*}
The above quantity $\mathbb{P}\left(\| V \|_{\infty} > \kappa \sigma \sqrt{2n\log (p/\eta)}\right)$ is smaller than $\eta$, if the parameter $\kappa$ is such that
\begin{equation*}
 p \exp \left( - \frac{1}{p}\sum_{j=1}^p \frac{ \kappa^2  n\log (p/\eta) }{\Omega_{j,j}}\right) = \eta,
\end{equation*}
or equivalently $\kappa = \sqrt{ \frac{1}{n}\frac{p}{\sum_{j=1}^p\Omega_{j,j}^{-1}  }}$. This is the announced result.
\end{proof}

\begin{proof}[Proof of Theorem~\ref{thmpreliminarylasso}]
We have, for any $\mu>0$,
\begin{equation*}
	\left\|A'A\left(\hat{\beta}^{DS}_{\mu} - \beta^{*}\right)\right\|_{\infty}
	\leq
	\left\|\left(A'A-X'X\right)\left(\hat{\beta}^{DS}_{\mu} - \beta^{*}\right)\right\|_{\infty}
        +
	\left\|X'X\left(\hat{\beta}^{LASSO}_{\mu}-\beta^{*}\right)\right\|_{\infty}.
\end{equation*}
Now, for the Dantzig Selector,
$$ \left\|X'X\left(\hat{\beta}^{DS}_{\mu}
             - \beta^{*}\right)\right\|_{\infty} \leq 2 \mu ,$$
with probability at least $1-\eta$,
provided that $\mu=2\sigma\sqrt{2n\log(p/\eta)}$.
Moreover,
$$
\|\hat{\beta}^{DS}_{\mu}
             - \beta^{*}\|_{1} \leq \|\hat{\beta}^{DS}_{\mu}\|_{1}
+ \|\beta^{*}\|_{1} \leq 2\|\beta^{*}\|_{1},
$$
implies that
$$ \left\|\left(A'A-X'X\right)\left(\hat{\beta}^{DS}_{\mu}
             - \beta^{*}\right)\right\|_{\infty}
\leq k\sigma \sqrt{2n\log(p)}.
$$
As a conclusion, with probability $1-\eta$,
\begin{equation*}
	\left\|A'A\left(\hat{\beta}^{DS}_{\mu} - \beta^{*}\right)\right\|_{\infty}
	\leq
	4 \sigma\sqrt{2n\log(p/\eta)} + k \sigma\sqrt{2n\log(p)}
	\leq
	(4+k)\sigma\sqrt{2n\log(p/\eta)}.
\end{equation*}
\end{proof}

\vspace{0.5cm}

\noindent {\bf Acknowledgment.} We would like to thank Professor Peter B\"uhlmann for insightful comments and also for providing us the motif scores dataset. We also would like to thank Professors Arnak Dalalyan, Alexander Tsybakov, Nicolas Vayatis, Katia Meziani and Joseph Salmon for useful comments.

\bibliographystyle{alpha}
\bibliography{GenLasDan}

\newcommand{\etalchar}[1]{$^{#1}$}
\begin{thebibliography}{MVdGB09}

\bibitem[AG03]{Amini03Trans}
M.~Amini and P.~Gallinari.
\newblock Semi-supervised learning with an explicit label-error model for
  misclassified data.
\newblock In {\em Proceedings of the 18th IJCAI}, pages 555--560. 2003.

\bibitem[AH08]{L1MOH}
P.~Alquier and M.~Hebiri.
\newblock Generalization of l1 constraint for high-dimensional regression
  problems.
\newblock Preprint Laboratoire de Probabilit\'es et Mod\`eles Al\'eatoires (n.
  1253), arXiv:0811.0072, 2008.

\bibitem[Aka73]{aic}
H.~Akaike.
\newblock Information theory and an extension of the maximum likelihood
  principle.
\newblock In B.~N. Petrov and F.~Csaki, editors, {\em 2nd International
  Symposium on Information Theory}, pages 267--281. Budapest: Akademia Kiado,
  1973.

\bibitem[Alq08]{CSEL}
P.~Alquier.
\newblock Lasso, iterative feature selection and the correlation selector:
  Oracle inequalities and numerical performances.
\newblock {\em Electron. J. Stat.}, pages 1129--1152, 2008.

\bibitem[AZ05]{Ando05Trans}
R.~K. Ando and T.~Zhang.
\newblock A framework for learning predictive structures from multiple tasks
  and unlabeled data.
\newblock {\em J. Mach. Learn. Res.}, 6:1817--1853 (electronic), 2005.

\bibitem[Bac08]{BachGpLasso}
F.~Bach.
\newblock Consistency of the group lasso and multiple kernel learning.
\newblock {\em J. Mach. Learn. Res.}, 9:1179--1225, 2008.

\bibitem[BBC{\etalchar{+}}05]{Balcan05Trans}
M.~Balcan, A.~Blum, P.~Choi, J.~Lafferty, B.~Pantano, M.~Rwebangira, and
  X.~Zhu.
\newblock Person identification in webcam images: an application of
  semi-supervised learning.
\newblock In {\em ICML Workshop on Learning with Partially Classified Training
  Data}. 2005.

\bibitem[BM98]{Blum98Trans}
A.~Blum and T.~Mitchell.
\newblock Combining labeled and unlabeled data with co-training.
\newblock In {\em Proceedings of the 11th Annual Conference on Computational
  Learning Theory}, pages 92--100. 1998.

\bibitem[BRT09]{Lasso3}
P.~Bickel, Y.~Ritov, and A.~Tsybakov.
\newblock Simultaneous analysis of lasso and {D}antzig selector.
\newblock {\em Ann. Statist.}, 37(4):1705--1732, 2009.

\bibitem[BT04]{BT04}
M.~Beer and S.~Tavazoie.
\newblock Predicting gene expression from sequence.
\newblock {\em Cell}, 117:185--198, 2004.

\bibitem[BTW07]{BTWAggSOI}
F.~Bunea, A.~Tsybakov, and M.~Wegkamp.
\newblock Aggregation for {G}aussian regression.
\newblock {\em Ann. Statist.}, 35(4):1674--1697, 2007.

\bibitem[Bun08]{Bunea_consist}
F.~Bunea.
\newblock {\em Consistent selection via the Lasso for high dimensional
  approximating regression models}, volume~3.
\newblock IMS Collections, 2008.

\bibitem[Cat07]{manuscrit}
O.~Catoni.
\newblock {\em PAC-Bayesian Supervised Classification (The Thermodynamics of
  Statistical Learning)}, volume~56 of {\em Lecture Notes-Monograph Series}.
\newblock IMS, 2007.

\bibitem[CH08]{ChriMo7GpLass}
C.~Chesneau and M.~Hebiri.
\newblock Some theoretical results on the grouped variables lasso.
\newblock {\em Mathematical Methods of Statistics}, 17(4):317--326, 2008.

\bibitem[CLLL03]{CLLL03}
E.~Conlon, X.~Liu, J.~Lieb, and J.~Liu.
\newblock Integrating regulatory motif discovery and genome-wide expression
  analysis.
\newblock In {\em Proceedings of the National Academy of Science}, number 100,
  pages 3339--3344. 2003.

\bibitem[CS99]{Collin99}
M.~Collins and Y.~Singer.
\newblock Unsupervised models for named entity classification.
\newblock In {\em Proc. Joint SIGDAT Conf. on Empirical Methods in Natural
  Language Processing and Very Large Corpora}, pages 100--110. 1999.

\bibitem[CSZ06]{semi-sup}
O.~Chapelle, B.~Sch\"olkopf, and A.~Zien.
\newblock {\em Semi-supervised learning}.
\newblock MIT Press, Cambridge, MA, 2006.

\bibitem[CT07]{Dantzig}
E.~Cand\`es and T.~Tao.
\newblock The dantzig selector: statistical estimation when $p$ is much larger
  than $n$.
\newblock {\em Ann. Statist.}, 35, 2007.

\bibitem[CZA05]{Chap05Trans}
O.~Chapelle, A.~Zien, and H.~Akaike.
\newblock Semi-supervised classification by low density separation.
\newblock In {\em Proceedings of the Tenth International Workshop on Artificial
  Intelligence and Statistics}, pages 57--64. 2005.

\bibitem[DT07]{ArnakTsyb}
A.~Dalalyan and A.B. Tsybakov.
\newblock Aggregation by exponential weighting and sharp oracle inequalities.
\newblock {\em COLT 2007 Proceedings. Lecture Notes in Computer Science 4539
  Springer}, pages 97--111, 2007.

\bibitem[EHJT04]{Efron-LARS}
B.~Efron, T.~Hastie, I.~Johnstone, and R.~Tibshirani.
\newblock Least angle regression.
\newblock {\em Ann. Statist.}, 32(2):407--499, 2004.
\newblock With discussion, and a rejoinder by the authors.

\bibitem[FHHT07]{PCO}
J.~Friedman, T.~Hastie, H.~H\"ofling, and R.~Tibshirani.
\newblock Pathwise coordinate optimization.
\newblock {\em Ann. Appl. Statist.}, 1(2):302--332, 2007.

\bibitem[FL01]{FanLiScad}
J.~Fan and R.~Li.
\newblock Variable selection via nonconcave penalized likelihood and its oracle
  properties.
\newblock {\em J. Amer. Statist. Assoc.}, 96(456):1348--1360, 2001.

\bibitem[HCB08]{Barron2}
C.~Huang, G.~L.~H. Cheang, and A.~Barron.
\newblock Risk of penalized least squares, greedy selection and {L}$1$
  penalization for flexible function libraries.
\newblock Submitted to Ann. Statist., 2008.

\bibitem[Heb08]{Mo7SLasso}
M.~Hebiri.
\newblock Regularization with the smooth-lasso procedure.
\newblock Preprint LPMA, 2008.

\bibitem[Joa99]{Joach99TSVM}
T.~Joachims.
\newblock Transductive inference for text classification using support vector
  machines.
\newblock In {\em ICML}. 1999.

\bibitem[JRL09]{DASSO}
G.~James, P.~Radchenko, and J.~Lv.
\newblock Dasso: Connections between the dantzig selector and lasso.
\newblock {\em J. Roy. Statist. Soc. Ser. B}, 71:127--142, 2009.

\bibitem[KKL{\etalchar{+}}07]{interior}
S.~J. Kim, K.~Koh, M.~Lustig, S.~Boyd, and D.~Gorinevsky.
\newblock An interior-point method for large-scale l1-regularized least
  squares.
\newblock {\em IEEE Journal of Selected Topics in Signal Processing},
  1(4):606--617, 2007.

\bibitem[Kol09a]{KoltchDant}
V.~Koltchinskii.
\newblock The {D}antzig selector and sparsity oracle inequalities.
\newblock {\em Bernoulli}, 15(3):799--828, 2009.

\bibitem[Kol09b]{Koltchl1plus}
V.~Koltchinskii.
\newblock Sparse recovery in convex hulls via entropy penalization.
\newblock {\em Ann. Statist.}, 37(3):1332--1359, 2009.

\bibitem[Lou08]{KarimNormSup}
K.~Lounici.
\newblock Sup-norm convergence rate and sign concentration property of {L}asso
  and {D}antzig estimators.
\newblock {\em Electron. J. Stat.}, 2:90--102, 2008.

\bibitem[MB06]{MeinshBulhmConsistLasso}
N.~Meinshausen and P.~B{\"u}hlmann.
\newblock High-dimensional graphs and variable selection with the lasso.
\newblock {\em Ann. Statist.}, 34(3):1436--1462, 2006.

\bibitem[MB07]{MB07}
Lukas Meier and Peter B{\"u}hlmann.
\newblock Smoothing {$l_1$}-penalized estimators for high-dimensional
  time-course data.
\newblock {\em Electron. J. Stat.}, 1:597--615, 2007.

\bibitem[MVdGB09]{VanGpLass}
L.~Meier, S.~Van~de Geer, and P.~B{\"u}hlmann.
\newblock High-dimensional additive modeling.
\newblock {\em Ann. Statist.}, 37(6B):3779--3821, 2009.

\bibitem[MY09]{MeinYuSelect}
N.~Meinshausen and B.~Yu.
\newblock Lasso-type recovery of sparse representations for high-dimensional
  data.
\newblock {\em Ann. Statist.}, 37(1):246--270, 2009.

\bibitem[NMTM99]{Nigam98ENTrans}
K.~Nigam, A.~McCallum, S.~Thrun, and T.~Mitchell.
\newblock Text classification from labeled and unlabeled documents using em.
\newblock In {\em Mach. Learn.}, pages 103--134, 1999.

\bibitem[Sch78]{bic}
G.~Schwarz.
\newblock Estimating the dimension of a model.
\newblock {\em Ann. Statist.}, 6:461--464, 1978.

\bibitem[Tib96]{Tibshirani-LASSO}
R.~Tibshirani.
\newblock Regression shrinkage and selection via the lasso.
\newblock {\em J. Roy. Statist. Soc. Ser. B}, 58(1):267--288, 1996.

\bibitem[Vap98a]{Vapnik}
V.~Vapnik.
\newblock {\em The Nature of Statistical Learning Theory}.
\newblock Springer-Verlag, 1998.

\bibitem[Vap98b]{Vapnik98}
V.~Vapnik.
\newblock {\em Statistical Learning Theory}.
\newblock Wiley, New York, 1998.

\bibitem[vdG08]{VandeGeerSparseLasso}
S.~van~de Geer.
\newblock High-dimensional generalized linear models and the lasso.
\newblock {\em Ann. Statist.}, 36(2):614--645, 2008.

\bibitem[vdGB09]{VandeGeerConditionLasso09}
S.~van~de Geer and P.~B{\"u}hlmann.
\newblock On the conditions used to prove oracle results for the lasso.
\newblock {\em Elect. Journ. Statist.}, 3:1360--1392, 2009.

\bibitem[Wai06]{WainSelection}
M.~Wainwright.
\newblock Sharp thresholds for noisy and high-dimensional recovery of sparsity
  using l1-constrained quadratic programming.
\newblock Technical report n. 709, Department of Statistics, UC Berkeley, 2006.

\bibitem[WSP07]{WangTransSVM07}
J~Wang, X.~Shen, and W.~Pan.
\newblock On transductive support vector machines.
\newblock In {\em Prediction and discovery}, volume 443 of {\em Contemp.
  Math.}, pages 7--19. Amer. Math. Soc., Providence, RI, 2007.

\bibitem[XGL03]{Zhu99Trans}
Zhu X., Z.~Ghahramani, and J.~Lafferty.
\newblock Semi-supervised learning using gaussian fields and harmonic
  functions.
\newblock In {\em ICML}. 2003.

\bibitem[XP05]{XiaoTrans05}
G.~Xiao and W.~Pan.
\newblock Gene function prediction by a combined analysis of gene expression
  data and protein-protein interaction data.
\newblock {\em J. Bioinformatics and Computat. Biol.}, 3(6):1371--1389, 2005.

\bibitem[YL07]{garrotte}
M.~Yuan and Y.~Lin.
\newblock On the non-negative garrotte estimator.
\newblock {\em J. Roy. Statist. Soc. Ser. B}, 69(2):143--161, 2007.

\bibitem[ZBL{\etalchar{+}}03]{Zhou99TransGraph}
D.~Zhou, O.~Bousquet, T.N. Lal, J.~Weston, and B.~Schoelk\"{o}pf.
\newblock Learning with local and global consistency.
\newblock In {\em NIPS 16}. MIT Press, 2003.

\bibitem[Zou06]{AdapLassoZou}
H.~Zou.
\newblock The adaptive lasso and its oracle properties.
\newblock {\em J. Amer. Statist. Assoc.}, 101(476):1418--1429, 2006.

\bibitem[ZY06]{BiYuConsistLasso}
P.~Zhao and B.~Yu.
\newblock On model selection consistency of {L}asso.
\newblock {\em J. Mach. Learn. Res.}, 7:2541--2563, 2006.

\end{thebibliography}

\end{document}